\documentclass[a4paper]{article}

\usepackage[T1]{fontenc}
\usepackage[utf8]{inputenc}
\usepackage{float}
\usepackage{graphicx}
\usepackage{subcaption} 

\usepackage[thinlines]{easytable}
\usepackage{hyperref}
\usepackage{wrapfig}
\usepackage{dingbat}
\usepackage{color}
\usepackage{refcount}
\usepackage[table]{xcolor}
\usepackage{hhline}
\usepackage[toc,page]{appendix}

\usepackage{amsmath,amssymb,amsthm}
\usepackage{mathtools}

\usepackage[english]{babel}

\usepackage{natbib}

\usepackage{xcolor}

\makeatletter
\@addtoreset{equation}{section}

\makeatother

\theoremstyle{Theorem}
\newtheorem{theo}{Theorem}
\theoremstyle{Lemma}

\newtheorem{rem}[theo]{Remark}
\newtheorem{assumption}[theo]{Assumption}

\theoremstyle{Definition}

\theoremstyle{Corollary}
\newtheorem{cor}[theo]{Corollary}
\newtheorem{example}[theo]{Example}

\DeclareFontFamily{U}{matha}{\hyphenchar\font45}
\DeclareFontShape{U}{matha}{m}{n}{
      <5> <6> <7> <8> <9> <10> gen * matha
      <10.95> matha10 <12> <14.4> <17.28> <20.74> <24.88> matha12
      }{}
\DeclareSymbolFont{matha}{U}{matha}{m}{n}
\DeclareFontSubstitution{U}{matha}{m}{n}

\DeclareFontFamily{U}{mathx}{\hyphenchar\font45}
\DeclareFontShape{U}{mathx}{m}{n}{
      <5> <6> <7> <8> <9> <10>
      <10.95> <12> <14.4> <17.28> <20.74> <24.88>
      mathx10
      }{}
\DeclareSymbolFont{mathx}{U}{mathx}{m}{n}
\DeclareFontSubstitution{U}{mathx}{m}{n}

\DeclareMathDelimiter{\vvvert}{0}{matha}{"7E}{mathx}{"17}

\newcommand{\E}{{\mathbb E}}

\newcommand{\N}{{\mathbb N}}
\newcommand{\Z}{{\mathbb Z}}

\newcommand{\R}{{\mathbb R}}

\newcommand{\Tr}{\textnormal{Tr}}

\newcommand{\C}{{\mathcal C}}

\newcommand{\HH}{{\mathcal H}}

\newcommand{\Mtr}{\textnormal{M}^{\rm tr}}

\newcommand{\be}{\begin{equation}} 
\newcommand{\ee}{\end{equation}}

\numberwithin{theo}{section}

\begin{document}

\title{Quantifying deviations from separability in  space-time functional processes }

\author{Holger Dette, Gauthier Dierickx, Tim Kutta}

\maketitle
\begin{abstract}

The estimation of covariance operators of spatio-temporal data 
is in many applications  only computationally feasible under simplifying assumptions, such  as 
separability of the covariance into strictly temporal and  
spatial factors. 
Powerful tests for  this  assumption  have  been proposed 
in the literature.
 However, as real world systems, 
such as climate data are  notoriously inseparable, validating this assumption by statistical  tests, seems inherently questionable.  
In this paper we present an alternative approach: By virtue of 
separability measures, we quantify how strongly the data's
covariance operator diverges from  a  separable approximation. Confidence intervals  localize these measures with statistical guarantees. This method provides users with a 
flexible tool, to  weigh the computational gains of a separable model against the associated increase in bias.
As separable approximations we consider the established methods of partial traces and partial products, and develop  weak convergence 
principles for the corresponding estimators.  Moreover, we also prove such results for
  estimators of optimal, separable approximations, which are arguably of most interest in applications. In particular we
   present for the first time  statistical inference for this object, which has been  confined to estimation previously. 

Besides confidence intervals, our results encompass  tests for approximate separability. 
All methods proposed in this paper  are free of nuisance parameters and do neither  require
computationally expensive resampling procedures nor the estimation of  nuisance parameters. 
A simulation study underlines the advantages of our approach and its applicability is demonstrated by 
 the investigation of German annual temperature data.

\end{abstract}

\bigskip

%
%
%
%
%
%
%
%
%
%

\section{ Introduction }

Spatio-temporal data is ubiquitous  in many branches of modern 
statistics, such as medicine \citep{WorsMarNeelVanFristEvan1996,Lind2008,Skup2010},
 urban pollution \citep{krall2015},
climate research \citep{laurini2018,chattopadhyay2020},
 spectrograms derived from audio signals
  or geolocalized data  \citep{RabSchaf-book-1978,BelBarPetChed2011} 
 and perhaps 
most prominently in geostatistics \citep{MitGentGump2006,GneiGentGut2007}. For a review of spatio-temporal models with an extensive list of fields of applications see \cite{KyrJour-review-1999} and the references therein.
An important step in the statistical analysis  of this type of data is the estimation  of  the  covariance structure, 
which is required, for example, for  principal component analysis,
 prediction  or kriging.
Nowadays,  technological advances often  yield  ultra high-dimensional   space-time data  and
 a completely non-parametric estimation of the covariance  operator is computationally very expensive (if not impossible).
In some applications, such as the  European Union  earth observation program {\it Copernicus},
 where  $4.4$Pb data per year are collected, even the storage of the estimated covariance
 structure is not possible, for instance see \url{https://www.copernicus.eu/en/access-data}.

{A common approach to drastically reduce the number of parameters in the covariance matrix
and to obtain  computationally feasible procedures for large space-time data sets,
is the imposition of structural assumptions on the covariance  of the underlying process.}
 One important example in the analysis 
of space-time data is \textit{separability}, which means that 
the  covariance  
  is modeled as product of   a purely spatial and  a purely temporal part 
   \citep[see][]{KyrJour-review-1999,GneiGentGut2007,MarcGenton2007,Sher2011,CresWik-book-2011,Cres-book-2015}. 
    The simplifying power of this assumption is tempting in applications, even where 
separability is not physically self-evident.  As a consequence   statistical tests 
for the hypothesis of separability have been developed  to validate  this assumption where
 possible and to  discard it where not.  For example,  
\cite{MatYaj2004,ScaMar2005,Fuentes2006,CruFer-CasGon-Man2010} proposed tests based on spectral methods for stationarity, axial symmetry and separability. In
     \citet{MitGentGump2005,MitGentGump2006,LuZim2005}  likelihood ratio tests under the assumption of a normal distribution were  investigated. 
  Finite dimensional non-parametric tests for separability were considered in \cite{LiGentSher2007}.
  More recent work on testing for separability  in a functional, non-parametric framework  can be found in   the papers of \cite{ConstKokReim2017,AstPigTav2017,ConstKokReim2018}
  or \cite{BagDette2017}.

On the other hand  many authors argue that, in applications it is often implausible to assume the complete absence of any space-time interactions in the covariance, which however is essential to separability 
\citep[see, for example,][]{KyrJour-review-1999,CresHuang1999,Gnei2002,HuaHsu2004,JunStein2004,Stein2005,LuZim2005,LiGentSher2007,TinHuy2013}.  
 This concern is strengthened  by the fact that in most empirical studies found in the literature, the hypothesis of  separability 
  is rejected \citep[see][]{MitGentGump2006,Fuentes2006,GneiGentGut2007,AstPigTav2017,ConstKokReim2018,BagDette2017}.
  Consequently,  testing the hypothesis of exact separability might be questionable from a decision theoretical  point of view 
  (as one already knows  that it cannot be  satisfied).  While conceptual critiques haunt "exact separability"  in theoretical discussions, we want to highlight a more pragmatic interpretation of this assumption, common in applications: Here separability means that a separable  operator provides a reasonable  approximation of the true covariance structure \citep[for instance, see][]{MarcGenton2007,AstPigTav2017}.

In the present work, we develop statistical methodology, which corresponds to this latter point of view. Instead of testing the hypothesis of "exact separability", we quantify how strongly a covariance operator deviates from separability. This approach is not only more realistic than previous tests, but also offers greater flexibility to users, who can judge in their own case, wether a deviation from separability is too  large  to justify a separable model, or wether computational gains outweigh the ensuing imprecision.

Our  analysis of separability is based on the normed  difference between the covariance  operator and a separable approximation, which vanishes 
if and only if  exact separability holds. We apply these  measures to the empirical covariance operator from the data, to quantify the deviation 
from  separability of the underlying covariance structure. Naturally the outcome depends critically on the particular method of separable approximation.
We examine three  types of approximations from a  general perspective. Two of them  have recently attracted attention in the problem of testing for exact separability, namely  partial traces \citep{ConstKokReim2017,ConstKokReim2018,AstPigTav2017}  and partial products \citep{BagDette2017}.  In order to bridge the aforementioned gap between earlier theory of "exact separability" and an inseparable reality we develop weak invariance principles
to encompass  models with inseparable covariance operators. 
Our most important contribution consists in the examination of  the third type, the {\it optimal separable approximation}  and the derivation of statistical inference for this vital object. While suboptimal approximations - as treated in  previous  literature on tests - are reasonable tools to test for exact separability, they are
prone to  overestimate  deviations from separability. On the other hand,  optimal approximations do not suffer from this imprecision. 
 Until now optimal, separable approximations have been confined to finite dimensional estimation \citep{MarcGenton2007}, and inference seemed out of reach due to their implicit definition. However in this paper we present a novel, stochastic linearization, which is employed to prove a weak invariance principle for the empirical, sequential estimator of the distance between the covariance operator and  an optimal approximation. Hence our theory provides tools to develop  statistical inference for the current estimation practice.

More precisely, the new  weak invariance principles  for all three approximation types are combined with self-normalization techniques to construct asymptotically pivotal estimates for the measures of deviation from separability. Therewith we construct  confidence intervals, such that deviations from separability can be quantified with statistical guarantees.
As a    further application  we develop  (asymptotically) valid statistical tests for the null hypothesis
    that {\it ``the distance between the covariance and a  separable approximation is small''}. 
This hypothesis of {\it ``approximate separability''}  sets our approach apart from the body of existing tests, which exclusively focus on the 
more theoretical concept of exact separability.

From a mathematical point of view our approach has several challenges. First - in contrast to tests  for exact separability -  the 
investigation  of asymptotic  properties cannot be conducted under the assumption of separability (the null hypothesis) 
and weak convergence 
of the estimators (and their sequential versions) has to be established for any underlying (not necessarily separable) covariance structure. Secondly, even if weak convergence of an estimator can be established
in the general case, the limiting 
distribution  depends on several nuisance parameters, which are extremely difficult to estimate. While in some cases the  bootstrap 
might be a solution to this problem, its computational cost is exorbitant for larger data sets. Moreover, for dependent data, as considered in this paper, 
 resampling procedures  and  methods based on estimation of the nuisance parameters prove labyrinthine, depending on further tuning parameters, such as the block length in a multiplier bootstrap  \citep[see][]{buecher2016} or a bandwidth in the long-run variance estimator \citep[see][Chapter 16]{HorKokBook2012}, which are in turn difficult to adjust properly. 
We sidestep such concerns by virtue of self-normalization, which yields asymptotically pivotal statistics. For this purpose we introduce a sequential process of estimates of the covariance operator and study its weak convergence  under an appropriate topology.
Thirdly,  in contrast to the related literature on non-parametric tests of separability, our weak convergence results are proved not for the space of  Hilbert--Schmidt, 
but for the more restrictive space of trace-class operators, which  does not enjoy a Hilbert structure. As a consequence establishing CLTs and invariance principles becomes technically more challenging. The benefit of this approach 
is, that it allows statistical analysis of all three separability measures (partial traces, partial products and optimal approximations) simultaneously, in a fully functional framework. 
One consequence of this analysis is the extention of existing work on  partial traces, which so far has been confined to finite dimensional spaces and was hence 
 reliant on  preliminary projections of the observed functions.  In contrast our methodology does not require  any preprocessing  and is thus easy to implement.

The remaining part of  this paper is organized as follows. In Section~\ref{sec2}  we recall some basic facts  about Hilbert spaces and their tensor products
and afterwards introduce measures for deviations from separability. Next we specify our framework in Section~\ref{sec3} and then state our  main theoretical results. In particular we
establish  the weak convergence of a sequential version of the empirical covariance operator in the 
vector space of all  trace-class-valued continuous functions.  These results are used 
to prove the weak convergence of appropriately standardized estimates of the measures of deviation from separability with a distribution free limit. As applications 
we derive asymptotically valid confidence intervals for these measures and valid statistical tests for the hypothesis of {\it approximate separability}.
Section~\ref{Sec_FiniSamp} is devoted to the finite sample properties of the proposed methods, both in simulations and by virtue of a data example. Finally, all proofs and technical lemmas are gathered in the Appendix. 

%
%
%
%
%
%
%
%
%
%

\section{Deviations from separability in Hilbert spaces}\label{sec2} 

In this section we develop the mathematical tools to analyze the separability of linear operators. In Section~\ref{Subsec_Op-in-HSp} we provide a minimal background on spectral theory, as well as tensor product Hilbert spaces. For more details we refer the interested reader to the monographs \cite{WeidmannBook1980} and \cite{GohGolKaas-book-1993,GohGolKaas-Book-2003}. 
In Section~\ref{Subsec_Meas-of-Sepa}, we introduce different measures to quantify  separability. Separability is assessed by the norm
of the   difference between a covariance  operator and its separable approximation.

%
%
%
%
%
%
%
%
%

\subsection{Operators on Hilbert spaces} \label{Subsec_Op-in-HSp}

For $i=1,2$ let $(\HH_i, \langle\cdot, \cdot  \rangle_{\mathcal{H}_i})$  denote generic, 
real, separable Hilbert spaces, equipped with inner products $ \langle\cdot, \cdot  \rangle_{\mathcal{H}_i}$.
The space of bounded, linear operators mapping from $\mathcal{H}_1$, to $\mathcal{H}_2$,
i.e. such operators $T$, which fulfill
$$
\vvvert T\vvvert_\mathcal{L}:=\sup_{\|x\|_{\mathcal{H}_1}\le 1}\|T(x)\|_{\mathcal{H}_2}<\infty,
$$
will be denoted throughout our discussion by $\mathcal{L}(\HH_1, \HH_2)$. 
Notice that all operators in this space are continuous.

 An important subspace of $\mathcal{L}(\HH_1, \HH_2)$ consists of all compact operators, i.e., all operators which map the unit ball of 
 $\HH_1$ to a compact set in $\HH_2$.
Any such operator $T: \HH_1 \to \HH_2$ can be diagonalized as follows
$$
T [ x ]= \sum_{i=1}^\infty \sigma_i \langle e_i, x \rangle_{\HH_1} f_i, \quad x \in \HH_1,
$$
where $\{e_i\}_{i \in \N}$ is an orthonormal basis of $\HH_1$ and $\{f_i\}_{i \in \N}$ of $\HH_2$. The \textit{singular values of $T$},
 $\{\sigma_i\}_{i \in \N}$, form a decreasing sequence of non-negative numbers, 
 converging to $0$.  The decay rate of the singular values, which may be seen as a measure of $T$'s regularity, is captured by the $p$-\textit{Schatten-norm}
	\begin{equation*} 
		\vvvert T \vvvert_p 
	:=
		\Big( \sum_{i=1}^\infty \sigma_i^p \Big)^{1/p} . 
	\end{equation*}
The corresponding Schatten class consists of all operators with finite norm
	\begin{equation*} 
		\mathcal{S}_p(\HH_1, \HH_2)
	:= 
		\{T \in \mathcal{L}(\HH_1, \HH_2):  \vvvert T \vvvert_p <\infty\}.
	\end{equation*}
It is well known that for any compact operator $T$, the inequality
$
\vvvert T \vvvert_p \le \vvvert T \vvvert_q
$
holds, for all $p>q\ge 1$. In particular this entails 
$\mathcal{S}_p(\HH_1, \HH_2) \supset \mathcal{S}_q(\HH_1, \HH_2)$.

In the following we will be particularly interested in the classes $\mathcal{S}_1(\HH_1, \HH_2)$, 
the so called \textit{trace-class operators} and $\mathcal{S}_2(\HH_1, \HH_2)$, the 
\textit{Hilbert--Schmidt operators}. The name "trace-class" is derived from the 
concept of "trace", as known from square-matrices. It  can be extended to an operator $T \in \mathcal{S}_1(\HH, \HH)$ acting on a separable Hilbert space $(\HH, \langle \cdot, \cdot \rangle_\mathcal{H})$ as 
\begin{equation*} 
		\Tr [ T ]
	:=
		\sum_{i=1}^\infty \langle T [ e_i ], e_i \rangle_\mathcal{H}, 
\end{equation*}
where  $\{e_i\}_{i \in \N}$, is an arbitrary orthonormal basis of  $\HH$.  If the operator $T$ is positive and symmetric, it holds that $\Tr(T)= \vvvert T \vvvert_1$.  More generally for an operator $T:\HH_1 \to \HH_2$ the trace may be used to express the trace-class and Hilbert--Schmidt norms as follows

$$  
\vvvert T \vvvert_1=\Tr\big[\sqrt{T T^*} \big] \quad \textnormal{and} 
\quad \vvvert T \vvvert_2=\Tr\big[T T^* \big].
$$
Here $T^*$ denotes the adjoint operator of $T$ and $\sqrt{TT^*}$ the canonical square root operator as defined for instance in \cite{HorKokBook2012}. 
We now turn to tensor products of Hilbert spaces. For the Hilbert spaces $\HH_1, \HH_2$, we  denote by $ \HH_1 \odot \HH_2$ their  algebraic
tensor product, on which a unique bilinear form, denoted by {$\langle \cdot , \cdot  \rangle$}, can  be defined such that
$$ 
 \langle h_1 \otimes h_2, h_1' \otimes h_2' \rangle := \langle h_1, h_1' \rangle_{\HH_1} 
 \, \langle  h_2, h_2' \rangle_{\HH_2}
$$
for any $h_1, h_1' \in \HH_1$ and $h_2, h_2' \in \HH_2$. Completing  $ \HH_1 \odot \HH_2$ 
with respect to the induced norm yields again a Hilbert space, referred to as  the 
\textit{tensor product Hilbert space}, which is denoted by 
$\HH_1 \otimes \HH_2$. For details on tensor product spaces we refer  to \cite{WeidmannBook1980}.

The tensor product structure can be extended to the bounded, linear maps acting on 
$ \HH_1 \otimes \HH_2$. Suppose that two operators $A \in \mathcal{L}(\HH_1, \HH_1)$ 
and $B \in \mathcal{L}(\HH_2, \HH_2)$ are given. Then we may define  the map 
	$$
		A \otimes B: \HH_1 \otimes \HH_2 \to \HH_1 \otimes \HH_2
	$$
 for all $h_1 \otimes h_2 \in \HH_1 \otimes \HH_2$ with 
$h_1 \in \HH_1$ and $h_2 \in \HH_2$ as
	$$
		A \otimes B [ h_1 \otimes h_2 ] := A [h_1] \otimes B [ h_2 ]. 
	$$
Requiring $A \otimes B$ to be linear extends the map uniquely to the whole 
space $ \HH_1 \otimes \HH_2$. Its norm $\vvvert A \otimes B \vvvert_\mathcal{L}$ is easily seen to equal 
$\vvvert A \vvvert_\mathcal{L} \vvvert B\vvvert_\mathcal{L} $.
Operators in $\mathcal{L}( \HH_1 \otimes \HH_2, \HH_1 \otimes \HH_2)$, which can 
be decomposed as a tensor product of operators in $\mathcal{L}(\HH_1, \HH_1)$ and 
$\mathcal{L}(\HH_2, \HH_2)$ are called \textit{separable}. 


%
%
%
%
%
%
%
%
%

\subsection{Measures of separability } \label{Subsec_Meas-of-Sepa}

Separability is an attractive property in applications of high dimensional and functional data, as it cuts memory space and boosts computational speed. While plain separability hardly exist in practice, there are cases where an operator is on the verge of being separable. To make this notion quantitative, different measures of separability have been proposed, most of which inspect the normed difference of an operator and a separable approximation. 
Here we consider three ways to construct such approximations, which have gained particular attention in the literature, namely approximations by partial traces \citep{ConstKokReim2017,AstPigTav2017}, by partial products \citep{BagDette2017} and by optimal approximations \citep{vanLoanPits1993,MarcGenton2007}.

%
%
%
%
%
%
%
%
%

\subsubsection{Partial traces} \label{Subsubsec_PartTr}

The concept of  \textit{partial traces} plays an important role in quantum mechanics, where it is used to decompose so-called density matrices in their marginals while preserving the "consistency of expectation" of observables of the subsystem.
For linear algebraic properties of partial traces, see for instance  \cite{Bhat2003} and the references therein for the physical motivation. Its extension to separable Hilbert spaces can be constructed as follows. Let $\HH := \HH_1 \otimes \HH_2$ be the tensor product Hilbert space of the separable Hilbert spaces $\HH_1$ and $\HH_2$. We then define the  
partial traces  
\begin{equation} \label{partial_trace}
		\Tr_1: \mathcal{S}_1(\HH, \HH) \to \mathcal{S}_1( \HH_2, \HH_2) 
	\quad 
		\textnormal{and} 
	\quad 
		\Tr_2: \mathcal{S}_1(\HH, \HH) \to \mathcal{S}_1( \HH_1, \HH_1 )
\end{equation}
for any trace-class operators  $A \in  \mathcal{S}_1(\HH_1, \HH_1)$ and $B \in  \mathcal{S}_1(\HH_2, \HH_2)$ as 
$$
		\Tr_1 [ A \otimes B]
	= 
		\Tr [ A ] B \quad \textnormal{and} \quad \Tr_2 [ A \otimes B ]
	=
	 	A \Tr[ B ].
$$
Requiring the partial traces $\Tr_1, \Tr_2$,  to be linear, yields unique 
extensions to the whole space. 
 We make a brief 
remark on the partial traces, which highlights their usefulness in the study 
of separability (see Proposition~C.1 in \cite{AstPigTav2017}, as well as  \cite{ConstKokReim2018}). 

\begin{rem} \label{Rem_RemPropPartialTraceOp}
Suppose $C \in \mathcal{S}_1(\HH, \HH)$ is a trace-class operator. Then it holds that:
\begin{itemize}
\item[1.] If $C$ is separable 
$$ 
		\Tr [ C ] C
	= 
	 \Tr_2 [ C ] \otimes \Tr_1 [ C ] .
$$
\item[2.] The partial traces are Lipschitz-continuous with respect to the $\vvvert \cdot \vvvert_{1 }$-norm. 
\end{itemize}
\end{rem}
    Note that  in the finite dimensional setting partial traces are also continuous w.r.t.\ the $\vvvert \cdot \vvvert_2$-norm, since all norms are equivalent. 
     This contrasts with the infinite dimensional setting. Indeed, consider a (separable) operator
     $C \otimes D$, where $C \in \mathcal{S}_1( \HH_1, \HH_1)$ and   $D \in \mathcal{S}_1( \HH_2, \HH_2 )$ are positive definite. It is easy to show that both
    equalities $ \vvvert C \otimes D \vvvert_2 =  \vvvert C \vvvert_2 \, \vvvert D \vvvert_2$ and
$
            \vvvert \Tr_1 [ C \otimes D ] \vvvert_2
    =
          \Tr [ C ] \, \vvvert  D \vvvert_2
   =
         \vvvert  C \vvvert_1 \, \vvvert  D \vvvert_2
$
    hold.
    However, the trace-norm is in infinite dimensions not equivalent to the Hilbert--Schmidt norm.
Hence the partial traces are not continuous with respect to the Hilbert--Schmidt norm. This seemingly minute technicality has important implications for the weak convergence results in Section~\ref{Subsec_WeakConvBSp} below, since the more restrictive space of trace-class operators is not a Hilbert space anymore.  
By virtue of the partial traces, we define the marginals of an operator $C$  as
\begin{equation} \label{C_1&C_2}
		C^{ ( 1 ) } 
	:=
		\Tr_2[C] /\Tr[C]  
	\quad 
		\textnormal{and} 
	\quad 
		C^{ ( 2 ) }
	:=
		\Tr_1[C] ,
\end{equation}
where we  have to assume that $\Tr[C] \neq 0$. For a covariance operator this means that the corresponding data is not deterministic. 
	Notice that if $C$ is separable, Remark~\ref{Rem_RemPropPartialTraceOp} guarantees the equality
$C=C^{ ( 1 ) }  \otimes C^{ ( 2 ) } $. 
If $C$ is not separable, the marginal product   
\begin{equation} \label{ltr}
C^{\rm tr}:=C^{ ( 1 ) }  \otimes C^{ ( 2 ) }
\end{equation}
 may still be regarded as a separable approximation of $C$. Hence one can assess separability by considering the measure
	\begin{align} \label{dev} 
		\Mtr[C] &:=\vvvert  C - C^{\rm tr}\vvvert_2^2.
	\end{align}
Large values of 	$\Mtr[C] $ indicate that $C^{\rm tr}$ is a poor approximation of $C$ and hence suggest only a 
small degree of separability.

%
%
%
%
%
%
%
%
%

\subsubsection{Partial products}  \label{Subsubsec_PartProd}

The definition of \textit{partial products} was introduced (without this particular name) in the work of \cite{BagDette2017}. Partial products provide an alternative to the partial traces in the study of separability. They are defined on the larger space of all Hilbert--Schmidt operators $\mathcal{S}_2(\HH, \HH)$ instead of the trace-class operators $\mathcal{S}_1(\HH, \HH)$.  The partial products are bilinear maps
\begin{equation*}  
		P_1: \mathcal{S}_2(\HH, \HH) \times \mathcal{S}_2(\HH_1, \HH_1) \to \mathcal{S}_2( \HH_2, \HH_2) 
	\quad 
		\textnormal{and} 
	\quad 
		P_2: \mathcal{S}_2(\HH, \HH) \times \mathcal{S}_2(\HH_2, \HH_2) \to \mathcal{S}_2( \HH_1, \HH_1) 
\end{equation*}
defined for a separable element $A \otimes B \in \mathcal{S}_2(\HH, \HH)$ and $C_1 \in \mathcal{S}_2(\HH_1, \HH_1 )$, $C_2 \in \mathcal{S}_2(\HH_2, \HH_2 )$  as
\begin{equation}
P_1(A \otimes B, C_1)
	= 
		 \langle A, C_1\rangle_{\mathcal{S}_2 (\HH_1, \HH_1) } B
	 \quad 
		\textnormal{and} 
	\quad 
		P_2(A \otimes B, C_2)
	=  
		\langle B, C_2\rangle_{\mathcal{S}_2 (\HH_2, \HH_2) } A. \label{properties_partial}
\end{equation}		
If $A$ and $B$ are trace-class, we can express the partial products in terms of the partial traces, e.g., when $B \in \mathcal{S}_1 ( \HH_2, \HH_2 )$, then denoting the adjoint of an operator $A$ by $A^*$,  the identity

	$$ 
		P_1(A \otimes B, C_1)
	=
		 \Tr_1\big[(A^* C_1) \otimes B \big]
	$$
	holds.
It was shown in \cite{BagDette2017} that  the operators $P_1$ and $P_2$ are well-defined, bi-linear and continuous with 
	\begin{align*}
& 
	\langle B, P_1(C,C_1) \rangle_{\mathcal{S}_2 (\HH_2, \HH_2)}
	= 
	\langle C,  C_1 \otimes B \rangle_{\mathcal{S}_2 (\HH, \HH)}
	\\
& 	
	\langle A, P_2(C,C_2) \rangle_{\mathcal{S}_2 (\HH_1, \HH_1)}
	=
	 \langle C,  A \otimes C_2 \rangle_{\mathcal{S}_2 (\HH, \HH)}.
	\end{align*}
In principle the partial products allow the decomposition of a separable operator $C$ as follows
$$C=\frac{P_2(C, \Delta_2)  \otimes P_1(C, \Delta_1)}{ \langle C, \Delta_1 \otimes \Delta_2 \rangle_{\mathcal{S}_2 (\HH, \HH)}},  $$
for any $\Delta_i \in \mathcal{S}_2(\HH_i, \HH_i)$ ($i=1,2$), as long as the denominator is not degenerate. When $C$ is not separable the right side may be regarded as a separable approximation of $C$. Its Hilbert--Schmidt distance to $C$  is minimized   for a fixed $\Delta_2$ by $\Delta_1= P_1(C, P_2(C, \Delta_2))/\vvvert P_2(C, \Delta_2) \vvvert_2^2$. Thus  \cite{BagDette2017} suggest the following measure of separability
\begin{align} \label{dev_dette} 
M^{\rm prod}[C] :=M^{\rm prod}[C, \Delta_2]:= \vvvert  C -C^{\rm prod}\vvvert_2^2,
\end{align}
where
\begin{equation} \label{lprod}
C^{\rm prod}:=\frac{P_2(C, \Delta_2)  \otimes P_1(C, P_2(C, \Delta_2))}{\vvvert P_2(C, \Delta_2) \vvvert_2^2}. 
\end{equation}

Obviously the quantification of $C$'s separability still depends on the specific choice of $\Delta_2$, which is an infinite dimensional parameter. A natural, ad hoc solution could be to choose $\Delta_2$ as a finite rank approximation of the identity. Then $P_2(C, \Delta_2) \approx \Tr_2[C]$, and the corresponding separability measure would be similar to that in \eqref{dev}.

%
%
%
%
%
%
%
%
%

\subsubsection{Optimal approximations}  \label{Subsubsec_OptAppr}

The evaluation of the above separability measures has important implications for users: A small measure suggests, that $C$ is well approximated by a separable operator, be it constructed via partial traces or partial products. However the interpretation for larger values is unsatisfactory. Is the covariance operator $C$ indeed far away from separability, or  is it just compared to an unsuitable approximation? To address this question it is  natural to consider the \textit{optimal, separable approximations} of $C$, that is an operator
\begin{align*} 
C^{\rm opt} \in  \textnormal{argmin} \Big\{ \vvvert C-\tilde C\vvvert_2^2: \tilde C= A\otimes B, A \in \mathcal{S}_2(\HH_1, \HH_1), B \in \mathcal{S}_2(\HH_2, \HH_2) \Big\}.
\end{align*}
In finite dimensions the existence of such approximations  has been established  already in the last century, see e.g.\ \cite{vanLoanPits1993}. These results were employed to study space-time data by \cite{MarcGenton2007}. The existence proof of optimal, separable approximations can be easily extended to arbitrary, separable Hilbert spaces. Notice that we may identify $ \mathcal{S}_2(\HH, \HH) \cong \HH \otimes \HH = \HH_1 \otimes \HH_2 \otimes \HH_1 \otimes \HH_2$ and $\mathcal{S}_2(\HH_1 \otimes \HH_1,  \HH_2 \otimes \HH_2) \cong \HH_1 \otimes \HH_1 \otimes \HH_2 \otimes \HH_2$. There exists a  canonical, isometric isomorphism of Hilbert spaces
	$$
		\Pi : \HH_1 \otimes \HH_2 \otimes \HH_1 \otimes \HH_2 \to \HH_1 \otimes \HH_1 \otimes \HH_2 \otimes \HH_2, 
	$$
which is uniquely determined by the images of simple tensors, namely  
	$$
		h_1 \otimes h_2 \otimes h_1' \otimes h_2' \mapsto h_1  \otimes h_1' \otimes h_2 \otimes h_2', 
	$$
where $h_1, h_1' \in \HH_1$ and $h_2, h_2' \in \HH_2$. For any separable map $A \otimes B \in \mathcal{S}_2(\HH, \HH)$, its image $\Pi[A \otimes B]$ is a rank one operator in 
$\mathcal{S}_2(\HH_1 \otimes \HH_1,  \HH_2 \otimes \HH_2)$. Conversely the preimage of any rank one operator in  $\mathcal{S}_2(\HH_1 \otimes \HH_1,  \HH_2 \otimes \HH_2)$ under $\Pi$ is separable. As $\Pi$ is in particular norm preserving, to find an optimal, separable approximation $C^{\rm opt}$ of $C$ is equivalent to determine the optimal rank-one approximation of $\Pi[C]$. We accomplish this by calculating the largest  singular value $\gamma_1$ of $\Pi[C]$ and  corresponding vectors $e_1 \in \HH_1 \otimes \HH_1$ and $f_1 \in \HH_2 \otimes \HH_2$. Then  we define the separable approximation
\begin{equation} \label{lopt}
		C^{\rm opt}
	:=
		  \Pi^{-1}[\gamma_1 \cdot e_1 \otimes f_1],
\end{equation}
which equals $C$ in case of separability.
The corresponding distance measure takes a particularly simple form
\begin{equation} \label{dev_opt}
M^{\rm opt}[C]:=\vvvert C- C^{\rm opt} \vvvert_2^2= \vvvert \Pi[C]- \Pi[C^{\rm opt}] \vvvert_2^2 = \vvvert C \vvvert_2^2 - \gamma_1^2.
\end{equation}
This measure is more easily interpretable than those constructed by partial traces and partial products, because it informs the user 
precisely about the possibilities and limits of separable approximations.
Notice that in general multiple, optimal approximations may exist.  However if $\gamma_1 > \gamma_2$ holds, which as we will see is a very mild condition, the optimal approximation is already unique.

\begin{rem}
{\rm We briefly comment on the favorable computational aspects of optimal approximations, in the context of covariance operators.
  In the subsequent statistical applications we estimate the separability of an unknown covariance operator $C$, by considering the measure $M^{\rm opt}[\hat C_n]$, where $\hat C_n$ is the empirical covariance operator, defined in \eqref{Eq_DefEmpCovOp} below. Calculating $M^{\rm opt}[\hat C_n]$, i.e., the Frobenius norm and the largest singular value of $\Pi[\hat C_n]$ 	can be efficiently implemented \citep{GolVanLoan2013,CarDeg2018}. 
Moreover it is possible to numerically determine these entities without saving the whole operator $\Pi[\hat C_n]$, but only the observed data 
 \citep{Tze2013}. 
Such storage concerns are material in some applications, as  indicated in the introduction.
We also note that often 
  one of the spaces, say $\HH_1$, is low dimensional (for details see Remark~\ref{spattem}). Under such circumstances $\hat \gamma_1$ is the largest  eigenvalue of the low-dimensional matrix $\Pi[\hat C_n]^* \Pi[\hat C_n]$, which can be easily calculated.   } 

\end{rem} $ $\\[-6ex]

\begin{example}  \label{hd21}
{\rm 
We conclude this section by considering a small example to illustrate  similarities and differences of the three separability measures in the finite dimensional case. Notice therefore that $\mathbb{R}^4 \cong \mathbb{R}^2 \otimes \mathbb{R}^2$. We consider the matrix 
$$C(q):= \left( \begin{array}{rrrr}
2 & 0 & 1 & q \\
0 & 2 & q & 1 \\
1 & q & 2 & q \\
q & 1 & q & 2 \\
\end{array} \right) \in \mathcal{S}_2(\mathbb{R}^2 \otimes \mathbb{R}^2,\mathbb{R}^2 \otimes \mathbb{R}^2) $$
where $q \in [0,1]$ is a separability parameter.

 \begin{figure}[h]
\begin{center}
    \includegraphics[width=0.35\textwidth]{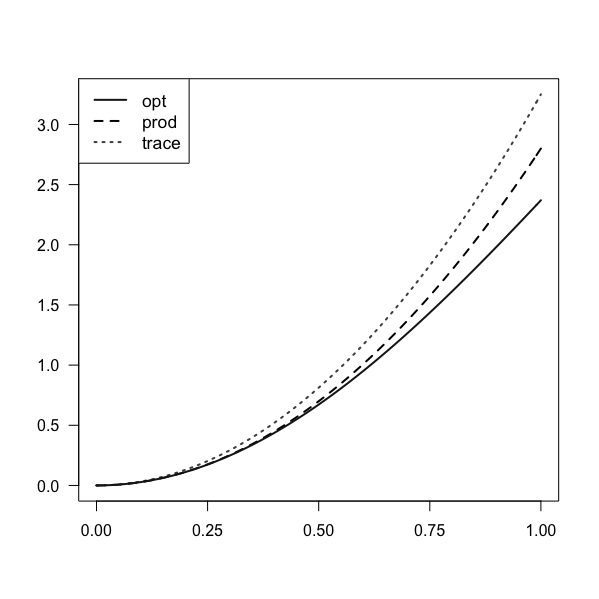}
  \caption{\textit{Comparison of the different measures $M^{\rm tr}, M^{\rm prod}, M^{\rm opt}$ in a matrix  example. For the partial products we have chosen $\Delta_2$ as the two dimensional identity matrix. \\[-8ex]} \label{Figure_1}
 } 
  \end{center}
\end{figure} 

The exact values of the measures \eqref{dev}, \eqref{dev_dette}  and \eqref{dev_opt} are calculated analytically  in Appendix~\ref{example}

 and given by 
 \begin{eqnarray*}
M^{\rm tr}[C(q)] &=& \vvvert C(q)-C(q)^{\rm tr} \vvvert _2^2=13/4 \cdot q^2  \\
M^{\rm prod}[C(q)] &=&  \vvvert C(q)-C(q)^{\rm prod} \vvvert _2^2=70/25 \cdot q^2 \\
M^{\rm opt}[C(q)]= &=&  \vvvert C(q) \vvvert _2^2-\gamma_1^2=10+3q^2 - \sqrt{9q^4+4q^2+100}
 \end{eqnarray*}
 
 For $q=0$, the matrix is separable and for larger $q$ the matrix gets increasingly inseparable. In Figure \ref{Figure_1} we observe a monotonic increase of all three measures in $q$. 
 Even for this  simple example the quality of the approximations $C^{\rm tr}, C^{\rm prod}, C^{\rm opt}$ varies considerably. In particular $M^{\rm opt}[C(1)]$ is only about $3/4$ of $M^{\rm tr}[C(1)]$ and for larger matrices this effect is generally much enhanced (see Section~\ref{Sec_FiniSamp}). 

Note that the  inequality $M^{\rm prod}[C(q)] \le M^{\rm tr}[C(q)]$
holds a priori for all $q$, since $C^{\rm tr}$ is the product of the marginals $C^{(1)} \otimes C^{(2)}$ as defined in \eqref{C_1&C_2}, while $M^{\rm prod}$ is the product of $C^{(1)}$ with an optimally adapted second factor.
}
\end{example}

%
%
%
%
%
%
%
%
%
%

\section{ Statistical Inference }\label{sec3}
{
In this section we proceed to the examination of separability in the covariance structure of space-time data. 
First we describe the statistical model considered in this paper and subsequently introduce the empirical covariance operator $\hat C_n$, as the natural, non-parametric estimator of the data's covariance operator $C$ (see equation \eqref{estimated_covariance_operator} for a precise definition).  Secondly,  we derive an invariance principle for the operator $\hat C_n$,
which is of interest in its own right.
  Since all measures $ M^{\rm x} [C]$  ($ \rm x \in \{tr, prod, opt \} $)
  introduced in Section \ref{sec2}  are functionals of $C$, corresponding estimates are obtained  replacing $C$ by  $\hat C_n$. Stochastic linearizations of the estimated measures $M^{\rm x}[\hat C_n]$ can be employed to prove their weak convergence to normal distributions. 
  However, the asymptotic variances depend in an intricate way on the unknown covariance operator $C$ and higher order moments 
  of the underlying model.
   To address this problem we use the concept of self-normalization to introduce standardized estimates of $M^{\rm x}[C]$, 
  that are asymptotically pivotal.
  These results have profound  statistical consequences, which are  explored in Section~\ref{sec33} and \ref{sec34}.
 }
 
%
%
%
%
%
%
%
%
%

\subsection{Space-Time Data}  \label{Subsec_SpaceTimeData}

 Let $\mu_1$, $\mu_{2}$ denote finite measures on the  interval $[0,1]$ and define 
 $H_1:=L^2_{\mu_{ 1 }}[0,1]$ and  $H_2:=L^2_{\mu_{ 2 }}[0,1]$   as  Hilbert spaces of square integrable functions with respect to the 
measure $\mu_{1} $ and $\mu_{2}$, respectively. The natural inner product on the  space  $H_{\ell}$ is defined as
$$
\langle f_1, f_2 \rangle_{\ell } := \int_0^1 f_1(t) f_2(t) \mu_{\ell }(dt)  \quad f_1,f_2 \in H_\ell~,~~\ell=1,2.
$$
For the ease of reading we will subsequently dispense with the indices of the inner products and norms, whenever there is no ambiguity. 
In the above situation, the tensor product Hilbert space $H:= H_1 \otimes H_2$ of 
spatio-temporal functions consists of all  square integrable functions  with respect to  the product measure $\mu_{1} \otimes \mu_{2}$, i.e. 
\begin{align} \label{spacegen} 
	H=\Big \{ h: [ 0 , 1]^{2}\to \R \, \Big | 
	\, \int_{[0,1]^{2}} h^2 ( s,t) (\mu_{1} \otimes \mu_{2} ) (ds ,dt) < \infty
	\Big\}.
\end{align} 
The corresponding inner product and norm on $H$ are given by  
$$
\langle h_1, h_2 \rangle =  \int_{[0,1]^{2}} h_{1}( s,t  ) h_{2}(  s,t) (\mu_{1} \otimes \mu_{2} ) (ds ,dt ) ~,~~ h_1, h_2  \in H 
$$
and $\| h \| = (\langle h, h \rangle )^{1/2}$, respectively.

\begin{rem} \label{spattem}

\rm{ In several applications of spatio-temporal statistics such as meteorology there exists a 
considerable asymmetry in the difficulty of acquiring spatial and temporal 
data. Whereas raising large numbers of densely timed 
observations at a certain location, say a weather station, 
is quite reasonable, establishing additional measuring stations across 
the country is often prohibitively expensive.  In this case  one can use 
 statistical models, which treat temporal  data as time-continuous functions
 (here continuity is achieved by  interpolations between the densely placed timepoints)
 while the spatial locations  are modeled by  a discrete parameter. 

 This situation corresponds to the choice of the Lebesgue measure for the measure $\mu_{2}$
 and  a discrete measure supported on $S$ points, say  $\{p_1,\ldots ,p_S\}$ for the measure $\mu_{1}$  in 
 \eqref{spacegen}. The inner product on the set  $H_1$ 
 is then given by  
 $$
\langle g_1, g_2 \rangle = \sum_{s=1}^S g_1(p_s) g_2(p_s) \quad g_1,g_2 \in H_1, 
$$
and the  tensor product Hilbert space $H$ in \eqref{spacegen} is given by 
	$$
	H=\Big \{ h: \{p_1,...,p_S\} \times [ 0 , 1 ] \to \R \, \Big |
	\, \sum_{ s =1}^S \int_0^1 h^2 ( p_s, t ) dt < \infty
\Big  \}.
	$$ 
}
\end{rem}

Suppose that  we observe a stretch $X_1, X_2, \ldots ,X_{n}$  of random 
functions from a  stationary,  $H$-valued  time series $(X_{i})_{i\in \Z}$, satisfying   $ \E \| X_1 \|^2 < \infty$. For the sake of a simple  presentation we assume that the 
random functions are centered, i.e. 
$$
\E X_i(s,t)=0, \quad    \forall s, t \in  \left[ 0 , 1 \right]^{2}.
$$
In practical applications this assumption is often unrealistic and it 
can be removed by empirically centering all the collected data without changing any of the following results
(see Remark~\ref{rem} below for more details). 
The second order structure of the random functions $X_i$ is captured by its covariance 
operator $C:= \E  X_1 \otimes X_1 \in \mathcal{S}_1(H,H) $ which  exists under the stated assumptions and is defined pointwise  by 
$$ 
		C ( h )
	:= 
		 \E X_1 \langle X_1, h \rangle \quad h \in H.
$$

Based on the sample  $X_{1}, \ldots , X_{n}$   the natural, non-parametric estimator of the covariance operator $C$ is given by
\begin{equation} \label{Eq_DefEmpCovOp}
		\hat{C}_n 	:= 	
		\frac{ 1 }{ n } \sum_{ i=1 }^{n }  \{X_i  \otimes   X_i\}.
\end{equation}
By virtue of $\hat C_n$ it is possible to estimate the separability of the true covariance operator $C$, as quantified by any of the measures $M^{\rm x}$, for  $\rm x \in \{tr, prod, opt\}$, simply by the plug-in $M^{\rm x}[\hat C_n]$. 
In the following  section we develop the analytic tools to examine the asymptotic behavior of these statistics.

%
%
%
%
%
%
%
%
%

\subsection{ Weak convergence in Banach spaces } \label{Subsec_WeakConvBSp}

We now devise the analytic tools for  the investigation of the statistics $ M^{\rm x}[\hat C_n]$. First we prove an invariance principle 
for a sequential version of the estimator $\hat C_n$, which is 
key for statistical self-normalization. We prove weak convergence on the space of trace-class-valued, continuous functions, i.e., sequential operators $L(\lambda)$, which live for some fixed $\lambda$ in $\mathcal{S}_1(H, H)$. 

	 This approach contrasts with the body of existing literature, where  weak convergence on the space of Hilbert--Schmidt operators $\mathcal{S}_2(H, H)$ is considered, see for instance the references in \cite{Mas2006} and more recently, \cite{ConstKokReim2018,BagDette2017}.
	  While  central limit theorems as well as invariance principles for this Hilbert space can be derived \citep{HoerKok2010,BerHorRic2013}, this approach is ineffective for the study of $M^{\rm tr}$, which has recently attracted much attention. The reason is, that $M^{\rm tr}$ is not well defined for all Hilbert--Schmidt operators (for details see Remark~\ref{Rem_GapProofTheoConstetal}). 
This problem can be sidestepped by projecting on a finite dimensional space, where all matrix norms are equivalent  \citep[see, for instance][for this approach]{AstPigTav2017}. Of course this requires the determination of an adequate projection space. For a fully functional approach we have to consider the covariance operator as a random element in $\mathcal{S}_1(H, H)$, which equipped with $ \vvvert  \cdot \vvvert_1$ is a (separable) Banach space \citep[see for instance][]{DunSchwartzBookPart1}. As pointed out by  \cite{Mas2006} this space might even be considered more natural for the study of weak convergence of empirical covariance  operators.\\

In order to specify the dependence structure of our data, we introduce the  notion of $\phi$-mixing,  as described, for  example in \cite{Samur1984}. To be precise for a  strictly
stationary sequence 
$(Y_j)_{j\in \Z}$ of random variables we  define 
	$$ 
		\mathcal{F}_{h,k }
	:=
		\sigma(Y_h , \ldots ,Y_k ),
	$$
	as the $\sigma$-algebra,  generated by $Y_h, \ldots  ,Y_k$  ($k \ge h$) and consider the {\it  $\phi$-mixing  dependence coefficients}
	\begin{align*}
		\phi(k) 
	:=
		\sup_{h \in \Z} \sup 
		\left\{ 
			{\mathbb{P}(F| E)} - \mathbb{P}(F)
			\, : \,  E \in \mathcal{F}_{1,h } , F \in \mathcal{ F }_{h+k, \infty}, \mathbb{P}(E)>0 
		\right\}.
	\end{align*}
	The sequence $(Y_j)_{j\in \Z}$  is called {\it  $\phi$-\textit{mixing}}, if $\lim_{k \to \infty }\phi(k) = 0$.  
	In order to prove  an invariance principle for a sequence
 of random elements  in the space $H$ defined in \eqref{spacegen} we make the following assumptions.


	

\begin{assumption}\label{assumption_1}~
\begin{itemize}
\item[(1)]   $(X_i)_{i \in \Z}$  is a strictly stationary time series  of centered, random functions in the  space  $H$ defined in \eqref{spacegen}. 
\item[(2)]  The sequence $(X_i)_{i \in \Z}$  is $\phi$-mixing, with $\phi(1)<1$ 
and such that 
$ \sum_{k=1}^\infty \sqrt{\phi(k)}<\infty. $
\item[(3)] There exists a basis $(e_q)_{q\in \N}$ of $H$ , such that 
$
\sum_{q=1}^{\infty} \sqrt[4]{\E | \langle X_1, e_q \rangle |^{4 }}  < \infty.
$
\item[(4)] The optimal separable approximation $C^{\rm opt}$ of $C$ is unique.
\end{itemize}
\end{assumption}

Notice that the third point in particular implies the existence of fourth moments of $X_1$, which follows by a simple calculation. 
The fourth point is much less restrictive than it might appear. In fact it is satisfied by a vast majority of geostatistical models considered in the literature. It is equivalent to assuming the strict inequality $\gamma_1 > \gamma_2$ for the two largest singular vaues of the restacked operator $\Pi[C]$ (see Section \ref{Subsubsec_OptAppr}). A simple, sufficient condition, that already entails it, is  the basic approximation requirement
$$\vvvert C - C^{\rm opt}\vvvert_2^2 < \frac{\vvvert C\vvvert_2^2}{2}.
$$

We now introduce a sequential estimator of the covariance operator $C$
	\begin{equation} \label{estimated_covariance_operator}
		\hat{C}_n ( \lambda )
	:= 	
		\frac{ 1 }{ n } 
		\sum_{ i=1 }^{ \lfloor n \lambda \rfloor }  ( X_i  \otimes   X_i )
		+\frac{n \lambda+1-\lceil n \lambda\rceil}{n} ( X_{\lceil n \lambda\rceil}\otimes X_{\lceil n \lambda\rceil} ),
	\end{equation}
	where $\lceil x \rceil = \min \{z \in Z \colon z \geq x\} $  and  
	 $\lambda \in [0,1]$. Note that $\lambda$
 determines the proportion of data used for estimation  and
 that  the estimate  in \eqref{Eq_DefEmpCovOp} is obtained as  $\hat{C}_n = \hat{C}_n (1)$. 
 The linear interpolation term on the right side of \eqref{estimated_covariance_operator} ensures that  $\hat{C}_n(\lambda)$ is continuous with respect to $\lambda$, which is important for technical reasons 
 (details are given in   Appendix~\ref{secA1}). A simple calculation shows that $\hat{C}_n$ is a trace-class operator 
 and by the discussion in Section~\ref{Subsec_Op-in-HSp}  a fortiori Hilbert--Schmidt. 
 	In  Theorem~\ref{Theo_funcCLT_EmpCovFunc} below, we provide an invariance 
	principle, which lies at the heart of the subsequent self-normalization approach. 

To state it in detail, we first recall the concept of a Gaussian measure and its associated Brownian motion on a Banach space.
Let $( F, \| \cdot \|_F )$ be a (separable) Banach space. Then, a measure $\mu_G$ on the Borel sets of $F$ is said to be centred Gaussian, if and only if $\mu_G \circ g^{ - 1 }$ is a centred Gaussian measure on $\R$ for all $g \in F^*$, the topological dual of $F$.  
Now, the vector space of $F$-valued, continuous functions denoted by
$$
\mathcal{C}(\left[ 0 , 1 \right], F):=\{f: \left[ 0 , 1 \right] \to F \, | \, f  \text{ continuous}\}~,
$$
 and equipped with the maximum norm 
\begin{equation} \label{maxnorm}
	\|f\|_\infty:=\sup_{\lambda \in \left[ 0 , 1 \right]} \|f(\lambda)\|_F,
\end{equation}
is again a (separable) Banach space. Let $\mu_G$ be a centred Gaussian measure on $F$ and let $G$ be the associated $F$-valued random variable. The  {\it Brownian motion $\mathbb{B}_G$ on $F$  corresponding to $G$}
is defined as  the  random function taking values in $\mathcal{C}(\left[ 0 , 1 \right], F)$ with the following properties:
\begin{itemize}
\item[(a)] $\mathbb{B}_G(\lambda)$ is distributed as $\sqrt{\lambda} G$ for any $\lambda  \in \left[ 0 , 1 \right]$.
\item[(b)] The increments of $\mathbb{B}_G$ are stationary, centred and independent. 
\end{itemize}
	An existence proof of the Brownian motion induced by a Gaussian measure on a (separable) Banach space can be found in \cite{Gross1970}.

 \begin{theo}\label{Theo_funcCLT_EmpCovFunc}

Let $(X_i)_{i\in \Z}$ be a sequence of random functions in $H$ satisfying 
Assumption~\ref{assumption_1}. Then there exists a
Gaussian random operator  $G \in \mathcal{S}_1 (H, H)$ such that 
$$
		\big \{ \sqrt{n } ( \hat{C}_n(\lambda)-C \lambda ) \big \}_{\lambda \in [0,1]}
	\stackrel { d } { \to }
			\big \{ \mathbb{B}_G(\lambda)\big \}_{\lambda \in [0,1]}
$$
in $\mathcal{C}(\left[ 0 , 1 \right], \mathcal{S}_1 (H, H))$. Here $\mathbb{B}_G$ is the Brownian motion on $\mathcal{S}_1 (H, H)$ corresponding to $G$. 
\end{theo}		

The proof of Theorem~\ref{Theo_funcCLT_EmpCovFunc} is deferred to the Appendix.  In the subsequent theorem we examine the Fr\'echet differentiability of the difference of a trace-class operator and its separable approximation.

\begin{theo} \label{Theo_FrechetDeriv}
Let $I \subset (0,1]$ be a closed interval and suppose that the covariance operator  $C$ satisfies Assumption \ref{assumption_1}, $(4)$. Then the maps
	\begin{itemize}
	\item[ i) ] $\mathbf{F}^{\rm tr}:   
\begin{cases}
   \mathcal{C}( I,  \mathcal{S}_1 ( H, H ) ) \to \mathcal{C}( I,  \mathcal{S}_1 ( H, H ) )\\
     L(\cdot ) \mapsto L  ( \cdot ) - L^{\rm tr}(\cdot),
		& 
\end{cases}$
	\item[ ii) ] $\mathbf{F}^{\rm prod}:\begin{cases}
    \mathcal{C}( I,  \mathcal{S}_2 ( H, H ) ) \to \mathcal{C}( I,  \mathcal{S}_2 ( H, H ) ) \\
     L(\cdot) \mapsto L  ( \cdot ) -L^{\rm prod}(\cdot),
		  & 
\end{cases}$
	\item[ iii) ] $\mathbf{F}^{\rm opt}:\begin{cases}
    \mathcal{C}( I,  \mathcal{S}_2 ( H, H ) ) \to \mathcal{C}( I,  \mathcal{S}_2 ( H, H ) ) \\
     L(\cdot) \mapsto  L  ( \cdot ) - L^{\rm opt}  ( \cdot)
	 ,  
		  & 
\end{cases}$
	\end{itemize}
 are Fr\'{e}chet differentiable in the function $\lambda \mapsto \lambda C $,
 where for fixed $\lambda \in I$ the operators $L^{\rm tr} (\lambda) $, $L^{\rm prod} (\lambda) $  and
	$L^{\rm opt} (\lambda) $ are defined in  \eqref{ltr}, \eqref{lprod} and \eqref{lopt}, respectively.
	 The derivatives are given in Section \ref{secA2} in the Appendix.

\end{theo}

Theorem~\ref{Theo_FrechetDeriv} implies that any of the three approximation-maps presented in Section \ref{Subsec_Meas-of-Sepa} is Fr\'{e}chet-differentiable. Since the Hilbert--Schmidt norm is also differentiable, the functional-$\Delta$-method 
yields the following, weak convergence result for the estimated separability measures.

\begin{cor} \label{cor_main}
Under the assumptions of Theorem \ref{Theo_funcCLT_EmpCovFunc}, it holds that for any $\rm x \in \{tr, prod, opt\}$
$$\Big\{ \sqrt{n}\big( M^{\rm x}[ \hat C_n(\lambda)]- M^{\rm x}[\lambda C] \big) \}_{\lambda \in I }\stackrel{d}{\to}  \big\{ \mathbb{B}(\lambda) \lambda \sigma^{\rm x}\big\}_{\lambda \in I }, $$
where $\mathbb{B}$ is a standard Brownian motion on the interval  $[0,1]$. The definition of the standard deviation $\sigma^{\rm x}$ is given in \eqref{sigmax}. 
\end{cor}

\begin{rem}\label{Rem_GapProofTheoConstetal}\rm 
	As mentioned in the preamble of this section, the statistic $\hat C_n - \hat C_n^{\rm tr}$ is  popular in the study of  separability. It has been investigated by  \cite{ConstKokReim2018} in the model described in Remark~\ref{spattem}.
	In their Theorem~2  it is indicated that  the limit of the operators $\sqrt{n}(\hat C_n - \hat C_n^{\rm tr})$ is Gaussian in the space $\mathcal{S}_2( H, H )$. Using their theorem, one could 
	 - in principle - derive the weak convergence of the corresponding measure $M^{\rm tr}[\hat C_n] $ 
	 via a continuous mapping argument, but the proof of this result is not complete. 
	More precisely, the partial traces are only continuous  with respect to the
Hilbert--Schmidt topology in finite dimensions, see  the discussion following Remark~\ref{Rem_RemPropPartialTraceOp}. 
However, in the functional setting continuity in an infinite dimensional space is required and crucial for an application of the functional-$\Delta$-method.
Thus, strictly speaking Theorem~2 in  \cite{ConstKokReim2018} is only correct for finite dimensional spaces.
However, Theorem~ \ref{Theo_funcCLT_EmpCovFunc} and  \ref{Theo_FrechetDeriv} and their  proofs in Section \ref{secA1} and \ref{secA2} provide  a solution to this problem. They establish weak convergence of $\sqrt{n}(\hat C_n - \hat C_n^{\rm tr})$ in $ \mathcal{S}_1 (H, H) $ 
as well as   Fr\'{e}chet differentiability of the map $L \mapsto L- L^{\rm tr}$ in $ \mathcal{S}_1 (H, H) $.
\end{rem}

\begin{rem}  \label{rem1} 
{\rm 
 There exists a large amount of literature  on CLTs  (often as a consequence 
of strong approximations)  for Banach space valued random variables under mixing conditions 
\citep[see, for example,][]{KuePhil1980,DehPhil1982,Deh1983}). Strong approximations imply  weak invariance principles, i.e.,  weak convergence of the process
$\big\{ \sqrt{n } ( \hat{C}_n (\lambda)- C \lambda ) \big\}_{\lambda \in [0,1] }$.  This holds in particular  if  $(X_{j})_{j\in \Z}$ is an absolute regular sequence with  mixing coefficients  $\beta (k)$ satisfying  similar conditions as
stated in Assumption~\ref{assumption_1}. 

Another popular dependency concept is  that of 
$L^p$-$m$-approximability, which extends $m$-dependence. It is frequently  used for  Hilbert space valued time series
\citep{HoerKok2010,BerHorRic2013} and is therefore potentially suited to the study of the separability measures $M^{ \rm prod }$ or  $M^{ \rm opt  }$, which are defined on the Hilbert space $\mathcal{S}_2(H,H)$.  However the measure $M^{ \rm tr }$ requires a Banach space setting, which
motivates the consideration of $\phi$-mixing in Assumption \ref{assumption_1}.   
We expect that  a CLT and perhaps an invariance principle for  
$\big \{ \sqrt{n} ( \hat{C}_n(\lambda)-C \lambda ) \big \}_{\lambda \in [0,1] }$ 
can be derived  in the case of $L^p$-$m$-approximability on $\mathcal{S}_1(H,H)$.
}
\end{rem}

%
%
%
%
%
%
%
%
%

\subsection{ Pivotal measures of the deviation from separability }  \label{sec33}

In this section we assess the separability of the data's covariance structure, by plug-in estimates of the separability measures $M^{\rm x}[C]$ $(\rm x \in \{tr, prod, opt\})$ introduced in Section~\ref{Subsec_Meas-of-Sepa}.  Note that 
Corollary \ref{cor_main} implies that  the statistic $M_n^{\rm x}[\hat C_n] - M^{\rm x}[C] $ satisfies a central limit theorem
\begin{align} \label{hd8} 
\sqrt{n} \big( M^{\rm x}[\hat C_n] - M^{\rm x}[C] \big) 	\stackrel { d } { \to } ~  {\cal N} (0, (\sigma^{\rm x})^{2}),
\end{align}
where the variance $(\sigma^{\rm x})^{2}$ is defined in equation \eqref{sigmax} of the Appendix.
The normal distribution in \eqref{hd8} is  non-degenerate if $\sigma^{\rm x} >0$, which in particular requires that $C$  is inseparable. 
Otherwise,  an application of Theorem~\ref{Theo_funcCLT_EmpCovFunc} and  the continuous mapping theorem, show that 
\begin{align*} 
n  	M_n^{\rm x}[\hat C_n] \stackrel { d } { \to } ~   \vvvert  G \vvvert _{2}^{2} ,
\end{align*} 
 where   $G$ is Gaussian random variable in $\mathcal{S}_{1} (H, H)$.
Note that 
$
  \vvvert  G \vvvert _{2}^{2} \stackrel { d } { =}   \sum_{j=1}^\infty \eta_j Z_j^2, 
$
where $Z_{1},Z_{2}, \ldots $,  are independent, standard normal random variables and $\eta_{1}, \eta_{2},  \ldots ,$ 
eigenvalues of a certain covariance operator $Q$.    However, the variance $(\sigma^{\rm x})^{2}$ and the values  $\eta_{1}, \eta_{2},  \ldots ,$ 
depend in an intricate way on the dependence structure of the time series $(X_{j})_{j\in \Z}$  and higher order moments of $X_{j}$.
As these quantities are  intensely difficult to estimate, the development of  statistical inference for measures of separability,  which 
is  directly  based on the asymptotic normality of $M^{\rm x}[\hat C_n]$ is an extremely  hard task.

To avoid the problem of  estimating nuisance parameters we develop a self-normalized statistic, which 
has an  easily accessible, asymptotic pivotal distribution. In particular its quantiles can be readily simulated.  For this purpose we define for $\rm  x \in \{tr, prod, opt\}$  the self-normalization term 
\begin{align} \label{Eq_EstimVariance-L2norm}
		&\hat V_{n}^{\rm x}
	:=
			\Big\{\int_I \big(M^{\rm x}[\hat C_n( \lambda )] - \lambda^2M^{\rm x}[\hat C_n)]  \big )^2 d\nu(\lambda)\Big\}^{1/2},
\end{align}
where $\nu$ is a probability measure on a  closed interval  $I \subset (0, 1]$.  
This particular  type of self-normalization  was also considered  by  \cite{DetteKokStan2018} 
in the context of change point detection 
 in the mean function of a functional time series.
The following result shows that
the self-normalized statistic converges to a pivot.

\begin{theo} 
\label{thmmain}
Let  $\rm x \in \{tr, prod, opt\}$. Suppose that Assumption~\ref{assumption_1} holds and that the variance $(\sigma^{\rm x})^2$ defined in \eqref{sigmax} is positive. Then the following, weak convergence holds
\begin{align} \label{hd3}
 	{M^{\rm x}[\hat C_n] - M^{\rm x}[C] \over \hat V_n^{\rm x}}  	\stackrel { d } { \to } ~ 
 &
 	W
	:= \frac{ \mathbb{B}(1) 
	}{
	 \big\{ \int_I \lambda^2 \big(\mathbb{B}(\lambda)-\lambda \mathbb{B}(1) \big)^2 d\nu(\lambda) \big\}^{1/2},
	 }
\end{align}
where $\mathbb{B}$ denotes a Brownian motion on the interval $[0,1]$.
\end{theo}

The quantiles of the limiting distribution in Theorem \ref{thmmain} can easily be simulated.  Exemplarily, we show in Table~\ref{Table_TableQuantPivotRV} the 
 quantiles of the distribution of  $W$  if the measure $\nu$ is a uniform distribution on the sets $\{1/20, \ldots, 19/20\} $ and  $\{1/30, \ldots, 29/30\} $. We also note that the distribution is symmetric, which implies
 that $q_{\alpha} =-q_{1-\alpha} $.
\begin{table}[h!]
\begin{center}
\begin{tabular}{|c |c c c|}
\hline
$K$ & $\alpha=0.01$ & $\alpha=0.05$ & $\alpha=0.1$ \\ [0.5ex]
\hline\hline
$20$& 16.479 & 9.895 & 7.097\\
\hline
$30 $& 16.248 & 9.925 & 7.149 \\
\hline
\end{tabular}\\[2ex]
\caption{\textit{$(1-\alpha)$-quantiles of the distribution of  the random variable $W$ defined in \eqref{hd3}, 
where $\nu$ is the uniform measure on the discrete set $\{l/K: l=1,\ldots , K-1\}$, for $K=20,30$.}}
\label{Table_TableQuantPivotRV}
\end{center}
\end{table}

\begin{rem}  \label{rem2} 
\rm ~
\begin{itemize}
\item[(1)] Theorem \ref{thmmain} allows statistical inference for  the important separability measure $M^{\rm opt}[C]$, which due to its implicit definition has not been conducted previously. In particular we can derive statistical tests and confidence intervals for the approximations presented in \cite{MarcGenton2007}.

\item[(2)] Other normalizing factors could be used as well. For example, defining
$$
\tilde V_n^{\rm x}
	:=
			\int_I \big| M^{\rm x}[\hat C_n( \lambda )] - \lambda^2M^{\rm x}[\hat C_n]   \big| d\nu(\lambda) 
$$
we obtain (if $\sigma^{\rm x} > 0$) the weak convergence
$$
 { M^{\rm x}[\hat C_n] - M^{\rm x}[C]  \over \tilde V_n^{\rm x}}  	\stackrel { d } { \to } ~  { \mathbb{B}(1) \over  \int_I \lambda \big | \mathbb{B}(\lambda)-\lambda \mathbb{B}(1) \big |  d\nu(\lambda) }.
$$
\end{itemize}

\end{rem} 

The separability measures considered so far quantify how good an operator $C$ can be approximated by a separable version in absolute terms. However these entities are only interpretable if the user has some preconceived idea of how large $\vvvert C \vvvert_2^2$ is. If this knowledge is not readily available, the most natural approach seems to be a comparison of $ M^{\rm x}[ \hat C_n]$ with the estimated Hilbert--Schmidt norm, i.e., $\vvvert \hat C_n \vvvert_2^2$. Accordingly we consider the relative separability measures
\begin{equation*}  
M_{\rm rel}^{\rm x}[C]:= \frac{M^{\rm x}[C]}{\vvvert C \vvvert_2^2},
\end{equation*}
its corresponding estimates $M_{\rm rel}^{\rm x}[\hat C_n]$
and the  normalization factors
\begin{align*}  
		&\hat V_{{\rm rel},n}^{\rm x}
	:=
			\Big\{ \int_I \big(\lambda ^2 M_{\rm rel}^{\rm x}[\hat C_n( \lambda )] - \lambda^2M_{\rm rel}^{\rm x}[\hat C_n)]  \big )^2 d\nu(\lambda) \Big\}^{1/2}.
\end{align*}
The measure $M_{\rm rel}^{\rm opt}$ has been proposed by \cite{MarcGenton2007} for space-time matrices, who called it "separability approximation error index".  It is easy to see that $M_{\rm rel}^{\rm opt
}[C]\in [0,1)$.  We can now prove an analogue to Theorem \ref{thmmain} for the relative separability measures.

\begin{theo} 
\label{thmmain_rel}
Let  $\rm x \in \{tr, prod, opt\}$. Under the Assumptions of Theorem \ref{thmmain} the following weak convergence holds
\begin{align} \label{hd4}
 {M_{\rm rel}^{\rm x}[\hat C_n] - M_{\rm rel}^{\rm x}[C] \over \hat V_{{\rm rel},n}^{\rm x}}  	\stackrel { d } { \to } ~ &W,
\end{align}
where the random variable $W$ is defined in \eqref{hd3}. 
\end{theo}

%
%
%
%
%
%
%
%
%

\subsection{ Statistical consequences }  \label{sec34}

In this section we discuss some statistical applications of the above, weak convergence results.
We begin with the construction of confidence regions
for the absolute and relative deviation measures. Throughout this paper we 
denote by  $q_{\alpha}$  the $\alpha$-quantile of the distribution of the random variable  $W$ defined in \eqref{hd3}.

\begin{theo} \label{thm1} Let $\rm x \in \{tr, prod, opt\}$. If Assumption \ref{assumption_1} is satisfied and $\sigma^{\rm x}>0$, which is defined in \eqref{sigmax}, then the intervals

\begin{align} \label{hd5} 
	\hat I_{n}^{\rm x} 
&= 
	\left[ M^{\rm x}[\hat C_n] + q_{\alpha/2} \hat V^{\rm x}_{n }, \, 
		M^{\rm x}[\hat C_n]  + q_{1-\alpha/2} \hat V^{\rm x}_{n} \right] 
\\
	\hat I_{ {\rm rel} , n}^{\rm x }  
&= 
	\left[ M_{\rm rel}^{\rm x}[\hat C_n] + q_{\alpha/2} \hat V^{\rm x}_{{\rm rel}, n }, \, 
		M_{\rm rel}^{\rm x}[\hat C_n] + q_{1-\alpha/2} \hat V^{\rm x}_{  {\rm rel}, n}  \right]  \label{hd6} 
\end{align}
are asymptotic $(1-\alpha)$-confidence intervals for the measures $M^{\rm x}[C], M^{\rm x}_{\rm rel}[C]$, respectively.
\end{theo}

By the duality of confidence intervals and tests the result of Theorem \ref{thm1} can be employed to construct separability tests. In contrast  to the currently available literature \citep[see][among others]{Fuentes2006,MarcGenton2007,ConstKokReim2017,AstPigTav2017,
ConstKokReim2018,BagDette2017} 
we do not  examine the problem of testing for exact separability, because as discussed   in our introduction,
separability is rarely  met in applications. Instead we want to determine whether a separable approximation $C^{\rm x}$ deviates from the true operator $C$ substantially or not.  Therefore we propose to test the hypotheses of a {\it relevant (relative) deviation from separability}, that is 
\begin{eqnarray} \label{hypothesis}
		H_0 
	: ~~
		M^{\rm x}[C]& =& \vvvert  C - C^{\rm x} \vvvert_2^2 \le \Delta
	\quad \textnormal{vs.} \quad
		H_1:			M^{\rm x}[C]  > \Delta \\
		\nonumber & \\ 
				H_0 
	: ~
		M^{\rm x}_{\rm rel}[C] &=& {\vvvert  C - C^{\rm x} \vvvert_2^2   \over \vvvert  C \vvvert_2^2  } \le \Delta
	\quad \textnormal{vs.} \quad
		H_1:			M^{\rm x}_{\rm rel}[C]  > \Delta_{\rm rel} 
		\label{hypothesisrel}
\end{eqnarray}
where   $\Delta$, $\Delta_{\rm rel} >0$   are given thresholds of relevance expressing the greatest deviation from separability, which is 
still considered as irrelevant in  the application at hand.   We propose to reject the null hypothesis in \eqref{hypothesis}, whenever
	\begin{equation} \label{test}
	 q_{\alpha} {\hat V_n^{\rm x}} +	M^{\rm x}[\hat C_n]>  \Delta		,
	\end{equation}
where $M^{\rm x}$ is defined in Section~\ref{Subsec_Op-in-HSp}, the normalizing factor  $\hat V_{n}^{\rm x}$ in  \eqref{Eq_EstimVariance-L2norm} and 
$q_{\alpha}$ is the $\alpha  $-quantile of the distribution of the random variable $W$ defined in  \eqref{hd3}.  
Similarly,  the null hypothesis in \eqref{hypothesisrel} is rejected, whenever
	\begin{equation} \label{testrel}
	q_{\alpha} {\hat V_{\rm rel,n }^{\rm x}}	+  M^{\rm x}_{\rm rel}[\hat C_n] >  	 \Delta	 ~.
	\end{equation}

These tests are motivated by the duality between confidence intervals and tests as, for example,  described in
\cite{aitchison1964}. To be precise, note that it follows by similar arguments as given in the proof of Theorem~\ref{thm1} that
 the interval
 $$
\tilde I_{n}^{\rm x} = [ q_{ \alpha} \hat V_n^{\rm x} +  M^{\rm x}[\hat C_n] ~ ,  \infty) 
 $$
 defines a $(1-\alpha)$ one-sided confidence interval for  the measure $M^{\rm x}[C] $. Now, by the duality between confidence interval and tests, an asymptotic 
 level $\alpha$-test for the hypotheses in \eqref{hypothesis} is obtained by rejecting the null hypothesis, whenever
 $$
 [0, \Delta] \cap \tilde I_{n}^{\rm x}  =   [0, \Delta] \cap \  [q_{ \alpha}  \hat V_n^{\rm x} + M^{\rm x}[\hat C_n] ~ ,  \infty)  = \emptyset , 
 $$ 
 which is  obviously equivalent to \eqref{test}.  As the same  arguments apply  to the measure $M^{\rm x}_{\rm rel}[C]$ these considerations lead to the following result.

\begin{theo} \label{self-norm_theorem}
Suppose that the assumptions of Theorem \ref{thm1} hold and  that $\Delta>0$ ($\Delta_{\rm rel}>0$). The test 
\eqref{test}, \eqref{testrel}   is  a consistent asymptotic-level-$\alpha$ test for the hypotheses 
\eqref{hypothesis}, 
\eqref{hypothesisrel}, respectively.  
\end{theo}

\begin{rem} \label{rem}~ \rm
\begin{itemize}
\item[(1)]  {\textrm  It follows  from the proof of Theorem \ref{thm1} that the probability of rejecting the  null hypothesis in \eqref{hypothesis} by the test \eqref{test}
converges to $0$ if  $M^{\rm x}[C] < \Delta$,  to $\alpha$ if  $M^{\rm x}[C]= \Delta$  and to $1$ if  $M^{\rm x}[C] > \Delta$.  A similar statement can be made for the hypothesis
 \eqref{hypothesisrel} and the test \eqref{testrel}
\item[(2)] The methodology can easily be extended to test  the hypotheses
\begin{eqnarray*} 
		H_0 
	: ~~
		M^{\rm x}[C]&>  &\Delta
	\quad \textnormal{vs.} \quad
		H_1:			M^{\rm x}[C]   \le  \Delta~.
\end{eqnarray*}
A consistent asymptotic level $\alpha$-test is obtained by rejecting the null hypothesis, whenever
\begin{equation*}  
		{M^{\rm x}[\hat C_n] } +q_{1-\alpha}{\hat V_n^{\rm x}}\leq   \Delta	 ,
	\end{equation*}
Note that this formulation of the testing problem allows to decide for an approximately separable covariance structure at a controlled type-I error.
Similarly, a consistent and asymptotic level $\alpha$-test for the hypotheses
$H_0 : M^{\rm x}_{\rm rel}[C] > \Delta_{\rm rel}  $ vs.  $H_1:	M^{\rm x}_{\rm rel}[C]   \le \Delta_{\rm rel} $ is obtained by rejecting the null hypothesis 
whenever
		$M^{\rm x}_{\rm rel}[\hat C_n] + q_{1-\alpha}{\hat V_{{\rm rel}, n}}	 \le	\Delta_{{\rm rel}}.$
The proofs  of these properties  are   omitted for the sake of brevity.
\item[(3)] 
A careful inspection of the proofs in the Appendix shows that all   results presented in this section  are also valid  without the assumption 
  $\mathbb{E}X_n(s,t) \equiv 0$. In this case one has to replace each  $X_n$ in the definition of $\hat C_n(\lambda)$ in \eqref{estimated_covariance_operator}, by the empirically 
 centred version $X_n-\bar X$, where $\bar{X}$ is the average function of  $X_1,\ldots  ,X_n$.
 }
\end{itemize} 

\end{rem}

%
%
%
%
%
%
%
%
%
%

\section{ Finite sample properties } \label{Sec_FiniSamp}
	
%
%
%
%
%
%
%
%
%

\subsection{ Simulations }

In this section we demonstrate the applicability of our approach by virtue of a simulation study. The measures $\mu_1$ and $\mu_2$ in \eqref{spacegen} are chosen as uniform distributions  on the sets $\{0,1/(S-1),2/(S-1), \dots,1\}$ and $\{0,1/(T-1),2/(T-1), \dots,1\}$, respectively.  Similarly as in \cite{ConstKokReim2018}, we generate synthetic data, according to the following moving average model
\begin{equation}
X_k(s,t):= \sum_{s'=1}^S \exp\big\{ -b^2 (s-s')^2 \big\} [e_k(t, s')+e_{k-1}(t, s')] \quad k=1, \dots,n,  \label{simulated_data}
\end{equation}
where the spatio-temporal, random functions $e_0, \dots,e_n$ are i.i.d.\ realizations of a centred Gaussian process $e \in H$, with  covariance function
\begin{align} \label{Genitings_sigma}
\sigma(s,s',t,t'):=&\frac{1}{(a|t-t'|+1)^{1/2}} \exp\Big( -\frac{b^2|s-s'|^2}{(a|t-t'|+1)^c}\Big)  + c \Big(1-\frac{|s^2-s'^2|}{2}-(t-t') \Big)_+,
\end{align}
and  we choose $a=10$ and $b=5$. 
The choice of the covariance structure \eqref{Genitings_sigma} is similarly motivated in \cite{CresHuang1999}, who investigate classes of separable and non-separable functions. Related versions have been used  for instance in \cite{GneiGentGut2007, ConstKokReim2018} and \cite{BagDette2017}. As we can see by inspecting $\sigma$, the covariance function is separable for $c=0$  and the corresponding covariance operator $C$ of the data inherits this separability. 

\begin{figure}[H] 
\centering
\includegraphics[width=1\linewidth]{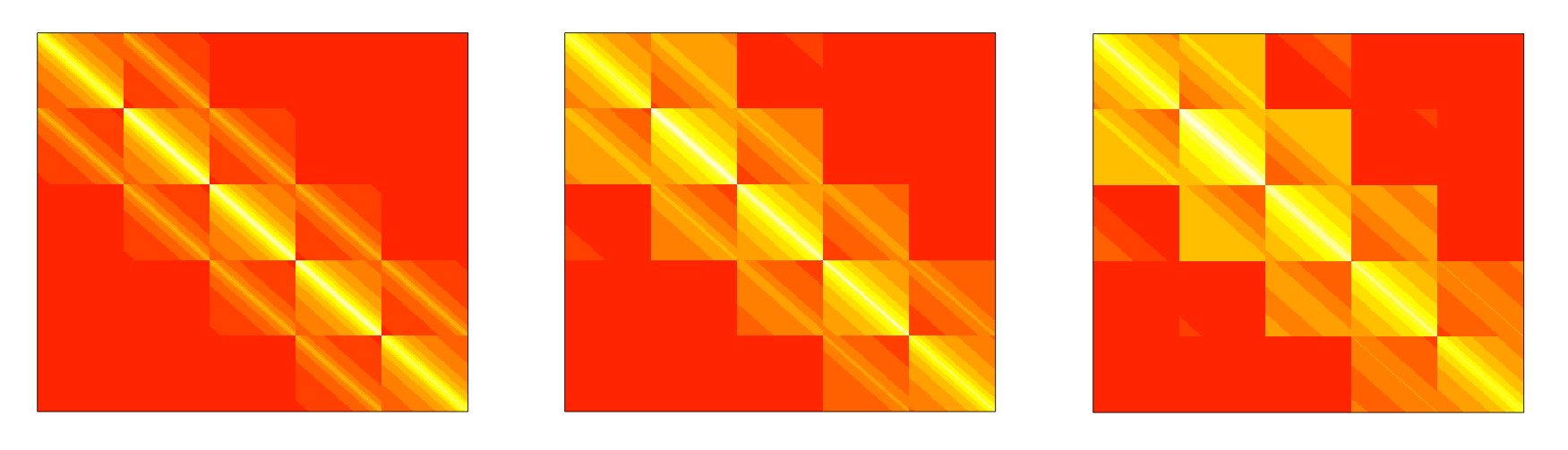}
\\[-2ex]
 \caption{\textit{The covariance operator $C$ in model \eqref{simulated_data} where
  $S=5$ ($25$ spatial blocks) and T=50, for different choices of the separability parameter $c$ (left $c=0$, middle $c=0.3$, right $c=0.6$).\label{different_c}}}
 \end{figure}

\begin{figure}[H]

\centering
\includegraphics[width=0.4\linewidth]{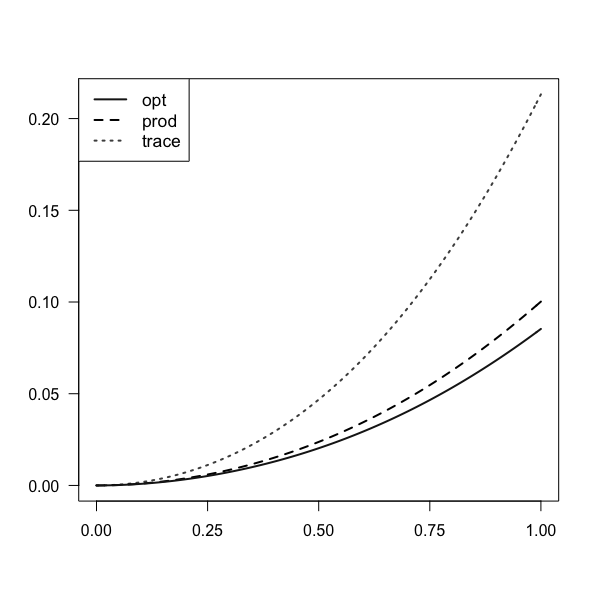}
 ~~~~~~~~~~~~
\includegraphics[width=0.4\linewidth]{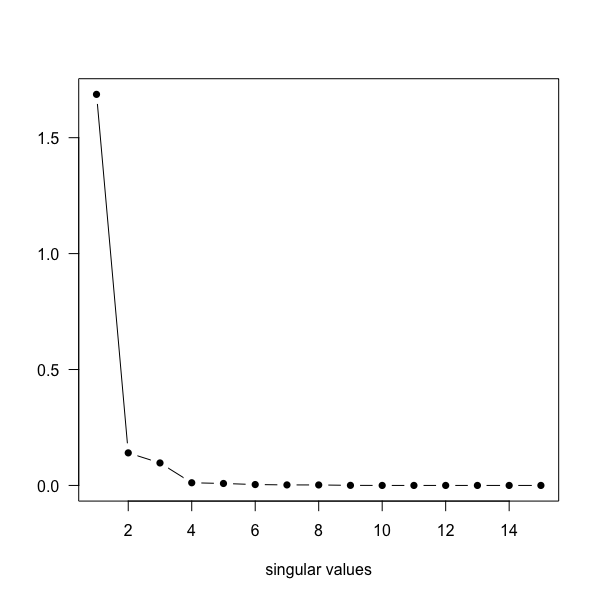}
\\[-2ex]
  \caption{\textit{ Left: Comparison of the separability measures $M^{\rm tr}, M^{\rm prod}, M^{\rm opt}$ for $S=5$  and $c \in [0,1]$ . \\
  Right: Singular values of $\Pi[C]$ for $S=5$ and fixed separability parameter $c=0.6$. \label{Figure_3 }}}
 \end{figure}

In Figure \ref{different_c} we can visually inspect the increasing deviation from separability: The left operator ($c=0$) is separable, i.e. all of the $25$ spacial blocks are scaled versions of one single time matrix. For $c=0.3$ and more strongly $c=0.6$ this homogeneity disappears.
The effect is most visible in the off-diagonal blocks which become asymmetric.  
In Figure \ref{Figure_3 } (left), we quantify the deviation from separability by virtue of our separability measures $M^{\rm x}[C]$. Notice that for the measure of partial products we have chosen the time matrix to be $\Tr_1[C]$ and optimize with respect to the space matrix (this corresponds to $\Delta_1= \operatorname{Id }$). In particular this implies that $M^{\rm tr}[C] \ge M^{\rm prod}[C]$.  Finally we can assess separability by inspection of the singular values of the 
restacked operator $\Pi[C]$. Recall that if $\Pi[C]$ is rank one (just one positive singular value) $C$ is separable. In Figure  \ref{Figure_3 } (right) we display the first $15$ singular values of $\Pi[C]$. The rapid decay after the first singular value indicates, that despite $C$'s observable inseparability, it is still  reasonably close to the optimal, separable approximation $C^{\rm opt}$. 
 
In the following we investigate the coverage probabilities of the asymptotic confidence intervals defined in \eqref{hd5}  and \eqref{hd6} for finite samples. For this purpose we simulate the data functions \eqref{simulated_data} for $T=50$ and  $S \in \{5,10, 15, 20\}$. All simulations are based on $5000$ simulation runs. The relative frequency of the events  
$M^{\rm x}[C] \in \hat{I}_n^{\rm x}$ and $M^{\rm x}_{\rm rel}[C] \in \hat{I}_{{\rm rel},n}^{\rm x}$ are displayed in Table \ref{Table_2}.

We observe a reasonable approximation of the coverage probabilities in most cases. 
First, the results are relatively stable with respect to the number  $S$  of spatial locations.
Second, the coverage probabilities  for the absolute differences  are  slightly too large and a better approximation is attained 
for the relative  differences (here in some case the coverage probabilities  are slightly too small). 
For the absolute measures, the confidence intervals for partial traces approximate  the prescribed $(1-\alpha)$-level
better than those for partial products and optimal approximations. In particular we frequently observe that $\mathbb{P} \{ M^{\rm opt}[C] \in \hat I_n^{\rm opt} \}  \ge \mathbb{P}\{ M^{\rm tr}[C] \in \hat I_n^{\rm tr} \}$. 
Despite this property  we almost invariably  observe that the interval $\hat I_n^{\rm opt}$ is is narrower than $ \hat I_n^{\rm tr}$ (these results are not displayed).
This is due to the superior quality of  optimal approximations compared  to those by partial traces.  This difference is visualized in the left panel of   Figure \ref{Figure_3 }, where we observe a sizable spread between the measures, already for  $c=0.6$. For larger $c$ the gap between $ M^{\rm tr}[C]$ and $ M^{\rm opt}[C]$ grows so rapidly, that they are hardly comparable at all. In contrast the measures $ M^{\rm opt}[C]$ and $ M^{\rm prod}[C]$ are very close in this example,  the corresponding intervals are similar and level-approximation is also comparable. However, it should be pointed out that there exist 
covariance operators  where the difference between  these measures is larger, as it has been indicated in  Example \ref{hd21} (see Figure \ref{Figure_1}). It will also become clear in the subsequent data example, that this difference potentially translates into noticably different confidence intervals.

\begin{table}[H]
{\footnotesize
\begin{center}
\begin{tabular}{|*{11}{c|}}
\hline  \multicolumn{3}{|c|}{ \cellcolor{gray!20} }  &
\multicolumn{4}{c|}{\cellcolor{gray!20} absolute } & \multicolumn{4}{c|}{\cellcolor{gray!20} relative } \\  \hline \cellcolor{gray!20}  & 
 $ \cellcolor{gray!20} n$ & \cellcolor{gray!20} $\alpha$& \cellcolor{gray!20} $S=5$ &  \cellcolor{gray!20} $S=10 $ & \cellcolor{gray!20} $S=15 $ & \cellcolor{gray!20} $S=20$  & \cellcolor{gray!20} $S=5$ & \cellcolor{gray!20} $S=10 $ &
\cellcolor{gray!20} $S=15 $ & \cellcolor{gray!20} $S=20$\\  \hhline{|=|=|=|=|=|=|=|=|=|=|=|}  \cellcolor{gray!20}  & \cellcolor{gray!20}
 $100$ & \cellcolor{gray!20} $5\%$ &  $0.969 $& $0.971 $& $0.971 $  &$0.971  $& $ 0.930 $
& $0.945 $ & $ 0.941$ & $0.941 $ \\ \hline  \cellcolor{gray!20}   & \cellcolor{gray!20}
   & \cellcolor{gray!20} $10\%$&$0.922 $  & $ 0.924$ & $ 0.926$& $0.924 $ & $ 0.858$ &$
0.876$ &$0.868 $ & $0.876 $  \\ \hline   \cellcolor{gray!20}  & \cellcolor{gray!20}  
   $200$ &  \cellcolor{gray!20} $5\%$ & $0.963 $ &$ 0.966$ &$ 0.967$&$0.963 $ &$ 0.955$&
$0.957 $& $0.953 $& $ 0.964 $ \\ \hline  tr \cellcolor{gray!20}   & \cellcolor{gray!20}  
   &  \cellcolor{gray!20} $10\%$& $0.909 $ & $0.917 $&$ 0.918$ & $  0.921$ & $0.888 $& $
0.899$ & $ 0.906$& $0.900  $ \\ \hline  \cellcolor{gray!20}   &  \cellcolor{gray!20} 
   $400$ &  \cellcolor{gray!20} $5\%$ & $0.957 $ &$0.958  $ &$0.960 $&$  0.959$ &$ 0.941  $&
$  0.961$& $0.964  $& $0.958   $ \\ \hline  \cellcolor{gray!20}   & \cellcolor{gray!20}  
   &  \cellcolor{gray!20} $10\%$& $ 0.899 $ & $0.905  $&$ 0.900 $ & $ 0.897  $ & $0.876  $& $0.902
 $ & $ 0.900 $& $0.904   $ \\  \hhline{|=|=|=|=|=|=|=|=|=|=|=|} \cellcolor{gray!20}  & \cellcolor{gray!20}  
$100$ &  \cellcolor{gray!20} $5\%$ & $0.979$ & $ 0.980$& $0.981 $  &$0.981  $& $0.932 $ &$
0.944$ & $0.944 $ & $ 0.940 $ \\ \hline  \cellcolor{gray!20}   & \cellcolor{gray!20}  
   &  \cellcolor{gray!20} $10\%$&$0.948$  & $ 0.934$ & $0.946 $& $0.944 $ & $0.856 $ &$ 0.870$
&$ 0.877 $ & $0.875 $  \\ \hline \cellcolor{gray!20}   &  \cellcolor{gray!20} 
   $200$ &  \cellcolor{gray!20} $5\%$ & $0.975 $ &$0.977 $ &$0.984 $&$0.985 $ &$0.955 $&
$0.938$& $ 0.932$& $0.931$ \\ \hline prod \cellcolor{gray!20}   & \cellcolor{gray!20}  
   &  \cellcolor{gray!20} $10\%$& $0.931 $ & $0.937 $&$0.947 $ & $0.952 $ & $0.894 $& $0.908 
$ &$ 0.919$& $0.917$ \\ \hline  \cellcolor{gray!20}   & \cellcolor{gray!20}  
   $400$ &  \cellcolor{gray!20} $5\%$ & $0.970  $ &$ 0.976 $ &$0.979 $&$ 0.982 $ &$ 0.940 $&
$ 0.962 $& $0.969  $& $0.979   $ \\ \hline  \cellcolor{gray!20}  &  \cellcolor{gray!20} 
   &  \cellcolor{gray!20} $10\%$&  $0.917  $ & $ 0.937 $&$ 0.940 $ & $0.942   $ & $0.875  $& $0.906
 $ & $ 0.917 $& $0.930   $ \\ 
 \hhline{|=|=|=|=|=|=|=|=|=|=|=|} \cellcolor{gray!20}   &  \cellcolor{gray!20} 
$100$ &  \cellcolor{gray!20} $5\%$ &$0.981 $  & $ 0.979$& $0.980$  &$0.981  $& $  0.928$ &$
0.940$ & $ 0.939$ & $0.938$ \\ \hline  \cellcolor{gray!20}  & \cellcolor{gray!20}  
   &  \cellcolor{gray!20} $10\%$&$0.939 $  & $ 0.937$ & $ 0.944$& $0.946 $ & $ 0.852$ &$
0.874$ &$0.861 $ & $0.871 $  \\ \hline  \cellcolor{gray!20}  & \cellcolor{gray!20}  
   $200$ &  \cellcolor{gray!20} $5\%$ & $0.975 $ &$0.966 $ &$0.985 $&$0.985 $ &$ 0.952$&
$0.959 $& $0.966 $& $0.966 $ \\ \hline opt \cellcolor{gray!20}  & \cellcolor{gray!20}  
   &  \cellcolor{gray!20} $10\%$& $0.932 $ & $0.937 $&$0.951 $ & $0.954 $ & $0.888 $ & $0.902
$ &$0.916 $& $0.916$ \\ \hline  \cellcolor{gray!20}  & \cellcolor{gray!20}  
   $400$ &  \cellcolor{gray!20} $5\%$ & $ 0.969 $ &$0.976  $ &$0.981 $&$ 0.984 $ &$0.938  $&
$ 0.959 $& $  0.966$& $ 0.979  $ \\ \hline  \cellcolor{gray!20}  & \cellcolor{gray!20}  
   &  \cellcolor{gray!20} $10\%$& $  0.918 $ & $ 0.935 $&$0.944  $ & $  0.945 $ & $ 0.870 $& $0.900
 $ & $ 0.915 $& $0.931   $ \\ \hline
\end{tabular}
\caption{\textit{ 
 Empirical coverage  probabilities of the confidence intervals  \eqref{hd5} and \eqref{hd6}. The data is generated according to model \eqref{simulated_data}, with separability parameter $c=0.6$. } \label{Table_2}}
\end{center}
}
\end{table}

%
%
%
%
%
%
%
%
%

\subsection{ Data Example }

In this section we employ the methodology developed in this paper to investigate the different measures of separability in a  data example. The data consists of daily temperature averages, published by the national meteorological agency ``Der  Deutsche Wetterdienst'' ( \url{ https://www.dwd.de/DE/Home/home_node.html} 
) at seven different stations in the cities of Bremen, Cottbus, Hohenpeissenberg, Karlsruhe, Magdeburg, Potsdam and  Schwerin. We consider data from the years $1893-1941$ and $1947-2008$, where a middle period is left out because of incomplete measurements during and immediately after WWII.

The data is smoothed over a Fourier basis with $41$ coefficients, which yields reasonable approximations of the temperature, while still reflecting general trends.
Furthermore it has been empirically centered and a linear trend of temperature increase has been removed to account for  climate change.

 \begin{figure}[H]
\centering
\includegraphics[width=0.4\linewidth]{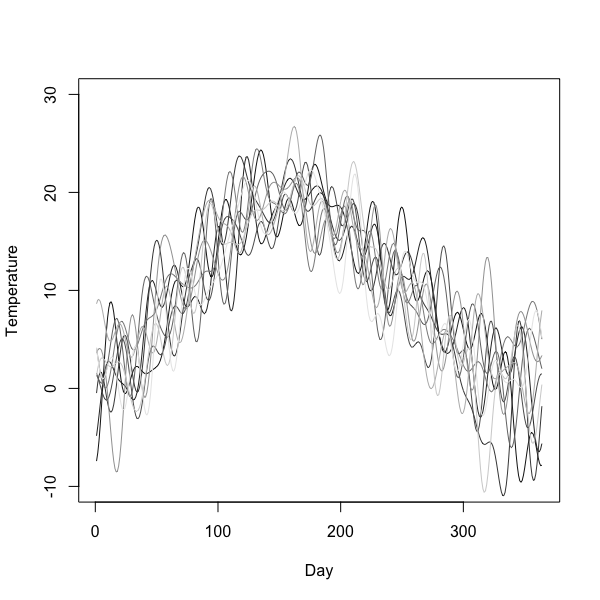}
 \label{fig:sub21b} 
~~~~ ~~~~ ~~~~ 
\includegraphics[width=0.4\linewidth]{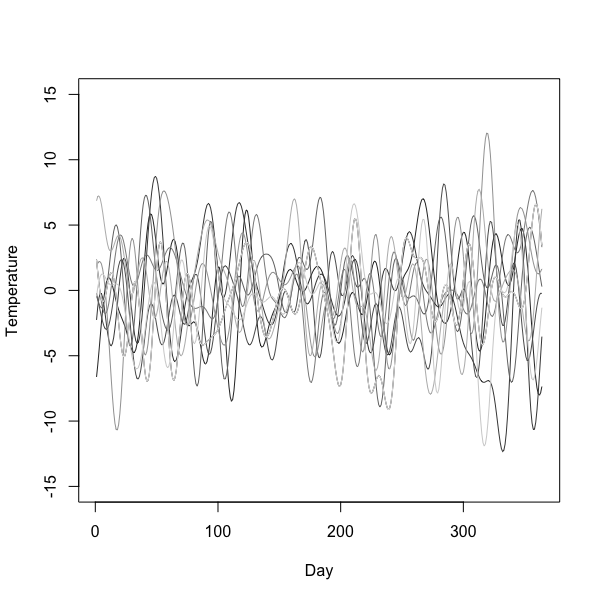}
         \caption{\textit{Temperature curves in the  first ten years  of observation at Hohenpeissenberg. Left: Observations. Right: Observations after 
        trend removal.} \label{Figure_4}}
 \end{figure}

\begin{figure}[H] 
\centering
\includegraphics[width=0.4\linewidth]{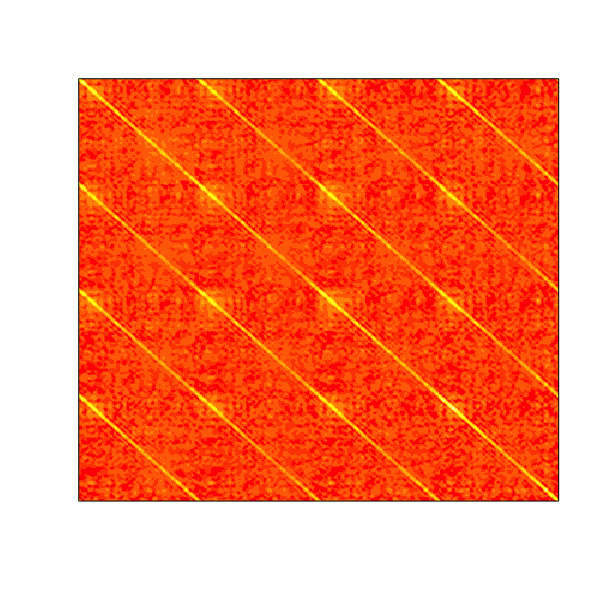}
~~~~~~~~~~
\includegraphics[width=0.4\linewidth]{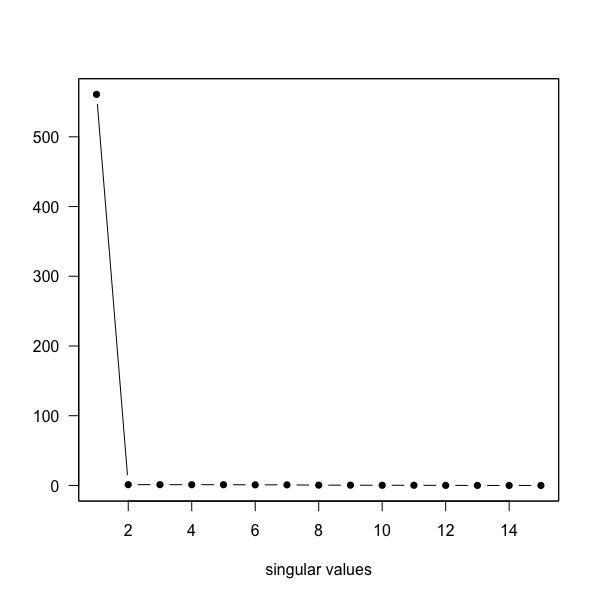}
 \\[-2ex]
         \caption{\textit{Left: Empirical covariance operator $\hat C_n$ for the cities Bremen, Hohenpeissenberg,  Karlsruhe and Potsdam. 
         Right: Singular values of $\Pi[\hat C_n]$. The fast decay after the first singular value indicates approximate separability. } \label{Figure_5}}
 \end{figure}

In Figure \ref{Figure_4} we illustrate  the effect of centering and  trend removal  at the first ten  observations at the location of Hohenpeissenberg (left before and right after detrending). We now turn to the investigation of the covariance operators of different collections of cities. Exemplarily we display in Figure \ref{Figure_5} (left) the covariance operator corresponding to the locations  Bremen, Hohenpeissenberg,  Karlsruhe and Potsdam. Compared to the covariance model considered in Section \ref{Sec_FiniSamp}, we observe a narrower concentration along the  diagonals of the block matrices. This almost-bandedness is also present in other geostatistical applications such as wind data  \citep[see, for example,  Figure 1 in][]{MarcGenton2007}. In Figure \ref{Figure_5} (right) we visualize the first fifteen singular values of the restacked operator $\Pi[\hat C_n]$, where the sharp decay after the first one already hints at  approximate separability. This impression will be reinforced by our subsequent investigations.

After these descriptive examinations we proceed to  statistical inference. We  consider different collections of cities and investigate the separability of the corresponding covariance operators by virtue of the empirical $95\%$-confidence intervals for the measures $M^{\rm x}[C]$ ($\rm x \in \{ tr, prod, opt\}$). In particular we examine the differences between approximation methods.

To make the analysis as comprehensible as possible, we restrict ourselves to 
the measures of relative deviation, which are more 
easily interpretable.  Moreover we do not display all possible combinations of cities, but confine ourselves to some illustrative examples.  In the following the $x$-axis represents  relative deviation in percent.

 \begin{figure}[H] 
 \captionsetup{format=hang}
\begin{subfigure}{0.5\textwidth}
\centering
\includegraphics[width=0.9\linewidth]{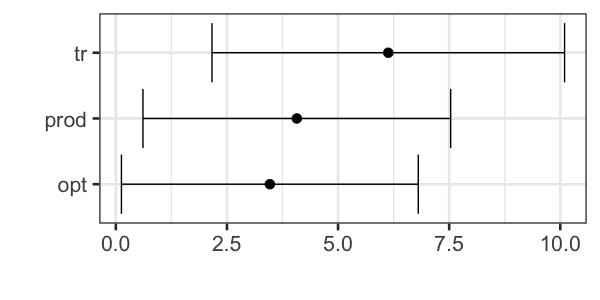}
 \\[-4ex]
  \caption{\textit{S=5:  Bremen, Hohenpeissenberg, Karlsruhe, Magdeburg, Potsdam }\label{fig:sub6a}}
 \end{subfigure}%
  \begin{subfigure}{0.5\textwidth}
 \centering
  \includegraphics[width=0.9\linewidth]{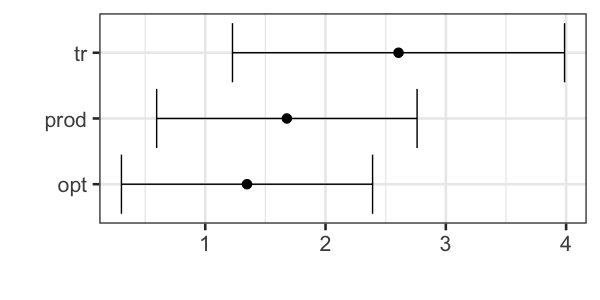}
\\[-4ex]
 \caption{\textit{S=4: Bremen,  Hohenpeissenberg, Karlsruhe, Potsdam} \label{fig:sub6b}}
 \end{subfigure}\\

\begin{subfigure}{0.5\textwidth}
\centering
\includegraphics[width=0.9\linewidth]{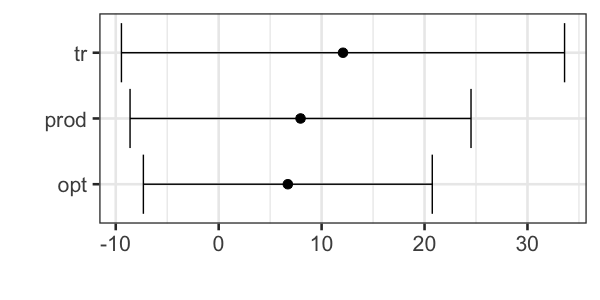}
\\[-4ex]
  \caption{\textit{S=6: Bremen, Hohenpeissenberg, Karlsruhe, Magdeburg,  Potsdam, Schwerin} \label{fig:sub6c}}
 \end{subfigure}%
  \begin{subfigure}{0.5\textwidth}
 \centering
  \includegraphics[width=0.9\linewidth]{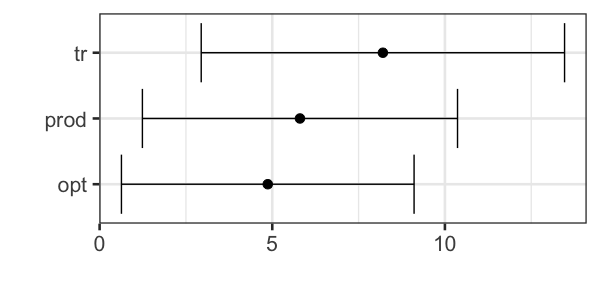}
\\[-4ex]
 \caption{\textit{S=4: Bremen, Karlsruhe, Potsdam, Schwerin, } \label{fig:sub6d}}
 \end{subfigure}\\

\begin{subfigure}{0.5\textwidth}
\centering
\includegraphics[width=0.9\linewidth]{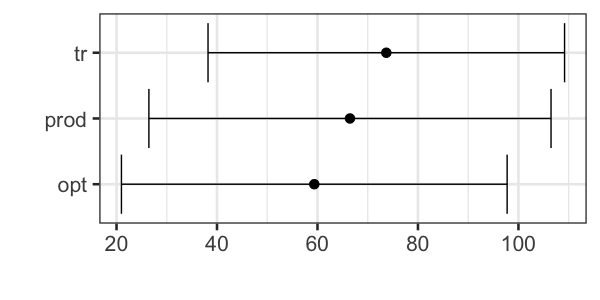}
 \\[-4ex]
  \caption{\textit{ S=7: Bremen, Cottbus ,  Hohenpeissenberg, Karlsruhe,
Magdeburg, Potsdam, Schwerin}\label{fig:sub6e}}
 \end{subfigure}%
  \begin{subfigure}{0.5\textwidth}
 \centering
  \includegraphics[width=0.9\linewidth]{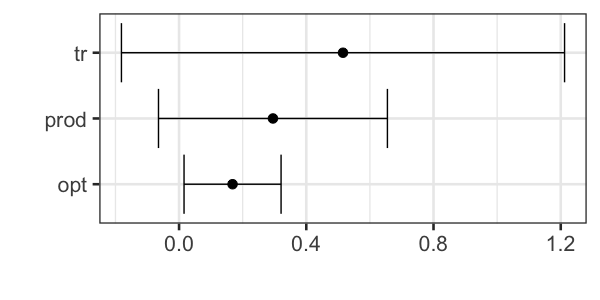}
 \\[-4ex]
 \caption{\textit{S=2: Hohenpeissenberg, Karlsruhe}\label{fig:sub6f}}
 \end{subfigure}\\
 \caption{ \textit{Empirical confidence intervals $\hat I_{{\rm rel},n}^{\rm x}$ for the relative separability  measures  $M^{\rm x}_{\rm rel}[C]$ $\rm ( x \in \{ tr, prod, opt\})$. 
 The $x$-axis represents  the  relative deviation in percent. } \label{Figure_6}}
 \end{figure}

In the top two panels of Figure \ref{Figure_6} (a) and (b), we observe deviation from separability but only to a small degree. In such cases (especially in (b)) it is reasonable to approximate the true covariance operator by a separable version, because the incurred imprecision is  small. The memory space, which is saved by separable approximations for $S=5$ amounts to $\approx 96\%$ and for $S=4$ to $\approx 94\%$, which is quite substantial. When we compare the measures of separability we notice considerable differences. For example in (b) we observe $M^{\rm tr}_{\rm rel}[\hat C_n] \notin 
 \hat I_{{\rm rel}, n}^{\rm opt} $, which underlines with statistical significance the differences of the measures. 

In the middle panels (c) and (d), we display examples of moderate deviations of separability. The confidence intervals for the relative measure of optimal approximations tend to be tighter than the other two (this is most visible in  (c), where $\hat I_{{\rm rel}, n}^{\rm opt} \subset \hat I_{{\rm rel}, n}^{\rm prod} \subset \hat I_{{\rm rel}, n}^{\rm tr}$). 

The bottom two examples  are  cases of either extreme  inseparability (e) and extreme separability (f). In (e) we observe that all of the relative measures assume values $>50\%$.  Here it is  not recommended  to approximate the covariance by a separable operator. Even linear combinations of separable operators  are not  suitable in this case. In contrast, on the  right we observe that a
separable approximation of the   covariance is reasonable. This may be useful despite the small number of spatial locations, as it still allows the  reduction of memory space to a quarter. 

\bigskip\bigskip

{\bf Acknowledgements}
This work has been supported in part by the
Collaborative Research Center ``Statistical modeling of nonlinear
dynamic processes'' (SFB 823, Teilprojekt A1, C1) of the German Research Foundation
(DFG). We are grateful to P. Kokoszka, M. Reimherr  and S. Tavakoli for pointing out important references and helpful discussions.

%
%
%
%
%
%
%
%

                                                                                                                                                                                                                                                                                                                                                                                                                \nocite{*}

\setlength{\bibsep}{2pt}

\newpage

\appendix

%
%
%
%
%
%
%
%
%

\section{Proofs and technical details}

%
%
%
%
%
%
%
%
%

\subsection{Details on the calculation of the separability measures} \label{example}

In this section we provide some details concerning the example of (in-)separability of the matrix 
$$C(q):= \left(\begin{array}{@{}cc|cc@{}}
    2 & 0 & 1 & q\\
    0 & 2 & q & 1 \\\hline
    1 & q & 2 & q \\
    q & 1  & q & 2 \\
  \end{array}\right) \in \mathcal{S}_2(\mathbb{R}^2 \otimes \mathbb{R}^2,\mathbb{R}^2 \otimes \mathbb{R}^2) $$
introduced in Section~\ref{sec2}. We will show in detail how each of the separability measures is derived, starting with the partial traces, where a  straightforward calculation yields
$$\Tr_1[C(q)]=\left(\begin{array}{rr}
2 & 0 \\
0 & 2   \\
\end{array} \right)+\left(\begin{array}{rr}
2 & q \\
q & 2   \\
\end{array} \right)=\left(\begin{array}{rr}
4& q  \\
q & 4  \\
\end{array} \right), \quad \Tr_2[C(q)]:= \left(\begin{array}{rr}
2+2 & 1+1  \\
1+1 & 2+2   \\
\end{array} \right)= \left(\begin{array}{rr}
4 & 2  \\
2 & 4   \\
\end{array} \right).$$
The trace of $C(q)$ is $8$.
Hence we have
$$C(q)^{\rm tr} =\Tr_2[C(q)] \otimes \Tr_1[C(q)]/\Tr[C(q)]=\left(\begin{array}{@{}cc|cc@{}}
    2 & q/2 & 1 & q/4\\
    q/2 & 2 & q/4 & 1 \\\hline
    1 & q/4 & 2 & q/2  \\
    q/4 & 1  & q/2 & 2  \\
  \end{array}\right)
 $$
It is now easy to see that $M^{\rm tr}[C(q)]=\vvvert C(q)-C(q)^{\rm tr}\vvvert _2^2=13/4 \cdot q^2$. For partial products and the parameter choice $\Delta_2 =  \operatorname{Id}_2$ (the two-dimensional identity matrix) we have, that $P_2 [ C(q), \Delta_2 ]=\Tr_2[C(q)]$, which is calculated above. Notice that $C(q)$ can be decomposed as the sum of the separable matrices
$$
C(q)=  \left(\begin{array}{rr}
1 & 0  \\
0 & 0   \\
\end{array} \right) \otimes  \left(\begin{array}{rr}
2 & 0  \\
0 & 2   \\
\end{array} \right)
+
\left(\begin{array}{rr}
0 & 1  \\
1 & 0   \\
\end{array} \right) \otimes  \left(\begin{array}{rr}
1 & q  \\
q & 1   \\
\end{array} \right)
+
\left(\begin{array}{rr}
0 & 0 \\
0 & 1   \\
\end{array} \right) \otimes  \left(\begin{array}{rr}
2 & q  \\
q & 2   \\
\end{array} \right).
$$

$P_1 [C, P_2(C, \Delta_2) ]$ can now be calculated using  bilinearity 
$$
	 P_1 [ C(q), \Tr_2[C(q)] ]=  4  \left(\begin{array}{rr}
2 & 0  \\
0 & 2   \\
\end{array} \right)
+
 4 \left(\begin{array}{rr}
1 & q  \\
q & 1   \\
\end{array} \right)
+
4 \left(\begin{array}{rr}
2 & q  \\
q & 2   \\
\end{array} \right)=\left(\begin{array}{rr}
20 & 8 q  \\
8 q & 20   \\
\end{array} \right),
$$
where we have used the defining property from \eqref{properties_partial} (the prefactors of the matrices are the Frobenius products of the space matrices and the partial trace $\Tr_2[C(q)]$). The squared Frobenius norm of $\Tr_2[C(q)]$ is $40$. Together this yields 
$$
		C^{\rm prod}(q) 
	= 
		\Tr_2[C(q)] \otimes P_1 [ C(q), \Tr_2[C(q)] ] / \vvvert \Tr_2 [ C(q) ] \vvvert_2^2 
	=
	 \left(\begin{array}{rrrr}
2\quad & 4/5 q & 1 \quad& 2/5 q  \\
4/5 q & 2  \quad& 2/5 q & 1 \quad \\
1\quad & 2/5 q & 2 \quad& 4/5 q  \\
2/5 q & 1  \quad & 4/5 q & 2 \quad \\
\end{array}\right)
$$
and hence a simple calculation reveals that $M^{\rm prod}[C(q)]= \vvvert C(q)-C(q)^{\rm prod}\vvvert _2^2=70/25 \cdot q^2$.
Finally we calculate the optimal approximation measure. First we notice that the squared Frobenius norm of $C(q)$ equals $ 20+6q^2$. Secondly, we calculate 
$$\Pi[C(q)] \Pi[C(q)]^*=   \left(\begin{array}{rrrr}
8 & 4 & 4& 8  \\
4 & 2  +2 q^2& 2  +2 q^2 & 4 +2 q^2 \\
4 & 2  +2 q^2 & 2  +2 q^2& 4 +2 q^2  \\
8 & 2  +2 q^2 & 4 +2 q^2 & 8  +2 q^2 \\
\end{array}\right),$$
the largest eigenvalue of which is given by $\gamma_1^2:= 3q^2 + \sqrt{9q^4+4q^2+100}+10$. As a consequence the measure of absolute deviation is given by $$M^{\rm opt}[C(q)]= \vvvert C(q)-C(q)^{\rm opt}\vvvert _2^2=\vvvert C(q)\vvvert _2^2-\gamma_1^2=10+3q^2 - \sqrt{9q^4+4q^2+100}.$$

%
%
%
%
%
%
%
%
%

\subsection{Proof of Theorem \ref{Theo_funcCLT_EmpCovFunc}}\label{Subsec_ProofTheofuncCLT} \label{secA1}

The proof of this theorem consists of two steps. First we establish 
a functional CLT, to demonstrate weak convergence of $\sqrt{n}(\hat C_n -C)$ in the space of trace-class operators. 
Subsequently by an application of Theorem~3.2 in \cite{Samur1987} we extend 
this result to an invariance principle. 
We note that the result of \cite{Samur1987} concerns weak convergence in the Skorohod space equipped with its Skorohod metric. Now, for continuous functions the Skorohod metric reduces to the supremum metric, see Chapter~3 of \cite{Bill-book-1999}. Since our random variables  as well as the limiting Brownian motion have  continuous sample paths, we obtain weak convergence in the space of continuous trace-class operators.
We begin by stating the CLT. 

\begin{theo}\label{Theo_CLT_EmpCovFunc}
	Let $(X_i)_{i\in \Z}$ be a sequence of random functions in $H$, defined in \eqref{spacegen},  satisfying Assumption~\ref{assumption_1}. Then there exists a Gaussian 
	random operator $G$ in $\mathcal{S}_1 (H, H) $ such that
	$$
		\sqrt{n } \big ( \hat{C}_n-C \big )
	\stackrel { d } { \to }
		G.
	$$\\[-8ex]
\end{theo}	

\noindent The proof of this theorem consists in the application of a CLT for 
Banach space-valued random variables, namely Theorem~4.4 of 
\cite{Samur1984}. In the following we are going to verify its four conditions:

\begin{itemize}
\item[(C1)] For any sequence of natural numbers $\{ r_n\}_{n \in \mathbb{N}}$ with
$r_n/n \to 0$ it follows that 
\begin{equation}\label{C1}
		\Big \vvvert\frac{1}{\sqrt{n}}\sum_{j=1}^{r_n} \big[ X_j \otimes X_j - C \big] \Big  \vvvert_{1} 	
	\stackrel { \mathbb{P} } { \to } 
		0.
\end{equation}
\item[(C2)] For any $\varepsilon>0$
\begin{equation}
	r_n \mathbb{P}\left\{
		 \vvvert X_1\otimes X_1-C\vvvert_1/\sqrt{n}>\varepsilon \right\}
 \to
 	 0. \nonumber
\end{equation}
\item[(C3)] For all $L \in \mathcal{L}(H,H)$ the following limit exists
\begin{equation*}
		\Phi(L)
	:= 
		\lim_{n \to \infty}  \E  \, \Tr \Big  (  \sum_{j=1}^{n} \left[ X_j \otimes X_j - C \right] \circ L\Big  )^2. 
\end{equation*}
\item[(C4)] There exists an increasing sequence of finite dimensional 
subspaces $F_k \subset \mathcal{S}_1(H, H)$, such that 
$$
	\lim_{k \to \infty} \sup_n  \E  \inf 
	\bigg  \{
		\Big  \vvvert \sum_{j=1}^{n}  \big( X_j \otimes X_j - C \big  ) -f \Big \vvvert_1 
		\, : \, 
		f \in F_k
	\bigg  \}
=
	0,
$$
where $\bigcup_{k \in \mathbb{N}} F_k $ is dense in $\mathcal{S}_1(H,H)$. 
\end{itemize}
We will now proceed to verify each of the claims.\\

In order to show (C1), we prove convergence to $0$ of the left side of 
\eqref{C1}  in expectation, which entails convergence in probability. 
By  Lemma~32 on p.\ 1116 of \cite{DunSchwartzBookPart1},  we obtain for 
the expected trace norm 
\begin{align}\label{C1_2}
	 \E \Big  \vvvert \frac{ 1 }{ \sqrt{ n } } \sum_{j=1}^{r_n} \left( X_j \otimes X_j - C \right)\Big  \vvvert_{1} 
	\le  
		\E \sum_{q \in \N } \frac{ 1 }{ \sqrt{ n } } \Big  \| \sum_{j=1}^{r_n} Y_{ j, q } \Big  \|=:E_{n}	,
\end{align}
where we have defined the centered random variables
\begin{equation}  \label{Y_{j,q}}
		Y_{j,q}
	:=
		X_j \langle  X_j, e_q \rangle  - \E X_j \langle  X_j, e_q \rangle  \in H. 
\end{equation} 
An application of  Jensen's inequality to   the right hand  side of \eqref{C1_2} yields
\begin{align}	\label{C1_3a}
 E_{n} \leq 	
		\sqrt{ \frac{r_n}{ n } } \sum_{q \in \N } 
			\Big\{ \E \Big \langle  \sum_{j=1}^{r_n} Y_{j,q}, \sum_{k=1}^{r_n} Y_{k,q}\Big \rangle  \frac{1}{ r_n } \Big \}^{1/2}.
\end{align}
Recalling the stationarity of the sequence $Y_{j,q} $, $j, q \ge 1$, 
inherited by  $X_j$, $j \ge 1$, we can simplify this expression to
\begin{align}
	\sqrt{ \frac{ r_n }{ n } } \sum_{ q \in \N } 
	\Big \{ \E \sum_{ | h | < r_n } \langle  Y_{0,q},  Y_{|h|,q}\rangle  \left( 1-\frac{ | h | }{r_n } \right)\Big \}^{1/2} 
\le
	\sqrt{\frac{r_n}{n}}\sum_{q \in \N } 
	\Big \{ 2 \sum_{ h \ge 0 } \left| \E \langle  Y_{0,q},  Y_{h,q} \rangle \right| \Big \}^{1/2}  
	\label{C1_3}
\end{align}
We will now exploit the $\phi$-mixing property of $Y_{h,q}$ to upper 
bound the covariance terms. We therefore use a special case of a 
 general result found in \cite{Dehl_Mik_Soren-Book-02}  on page~24, which yields 
\begin{equation}
	\left| \E  \langle  Y_{0,q},  Y_{h,q} \rangle \right| 
\le
	2 \sqrt{ \phi( | h | ) }  \E  \|Y_{ 0 , q } \|^2. \label{Dehling_Lemma}
\end{equation}
Plugging this bound in \eqref{C1_3} and observing \eqref{C1_3a} we obtain
\begin{align*}
E_{n } \le &
	\sqrt{  \frac{ r_n }{ n } } 
	\sum_{ q \in \N } 
	\Big \{ 2 \sum_{h \ge 0} 2 \sqrt{ \phi(h) }     \E \| Y_{ 0 , q } \|^2  \Big  \}^{ 1 / 2 }
= 	\sqrt{ \frac{ r_n }{ n } }  \sum_{ q \in \N }    \sqrt{  \E  \|Y_{ 0 , q } \|^2 } 
	\Big  \{ 4 \sum_{ h \ge 0 } \sqrt{ \phi( h ) }  \Big  \}^{1/2}.
\end{align*} 
Since by assumption $r_n /n \to 0$, it remains
to prove the finiteness of both sums on the right to establish (C1). The second one is finite 
due to our summability assumption on $\sqrt{\phi(h)}$. To see the 
boundedness of the first one we recall the definition of $Y_{0,q}$ 
in \eqref{Y_{j,q}}. By the binomial formula in Hilbert spaces we get
\begin{align*}
&\E \|Y_{0,q}\|^2 = \E \|X_0 \langle  X_0, e_q\rangle - \E X_0 \langle  X_0, e_q\rangle \|^2 
\le 2  \E \|X_0\|^2 \langle X_0, e_q \rangle ^2.
\end{align*}
 Applying the Cauchy--Schwarz inequality finally yields
$$ 
 \E \|X_0\|^2 \langle X_0, e_q \rangle ^2 \le \sqrt{\E \|X_0\|^4} \sqrt{\E\langle X_0, e_q \rangle ^4}.
$$
Noticing that by assumption 
$$\sum_{q \in \N} \sqrt[4]{\E | \langle X_0, e_q \rangle | ^4}< \infty,$$
 the proof of condition  (C1) is completed.\\
Verifying  (C2) is straightforward. Indeed, it follows immediately by an application of Markov's inequality and the dominated convergence theorem.\\
To establish condition (C3), let $L$ be an arbitrary bounded operator. By definition of the trace we observe that $\Phi(L)$ equals
\begin{align*}
& \lim_{n \to \infty} \frac{1}{n}  \E 
\Big \{\sum_{p\in \N}\langle [\sum_{i=1}^n X_i \otimes X_i -C] L [ e_q ] ,  e_p \rangle 
\sum_{q\in \N}\langle [\sum_{j=1}^n X_j \otimes X_j -C] L [ e_q ] ,  e_q \rangle\Big\}\\
=&\lim_{n \to \infty}  \frac{1}{n} \sum_{i,j=1}^n \E\Big\{\sum_{p\in \N}
\langle [ X_i \otimes X_i -C] L [ e_q ] ,  e_p \rangle \sum_{q\in \N}\langle [ X_j \otimes X_j -C] L [ e_q ] ,  e_q \rangle \Big\}.
\end{align*} 
We use the stationarity of the sequence $X_j $, $j \ge 1 $, to transform the last expression into
$$
		\lim_{n \to \infty}  
		 \sum_{ |h|<n} \Big ( 1-\frac{|h|}{n} \Big) 
		 \E \Big \{\sum_{p\in \N}
		 	\langle [ X_0 \otimes X_0 -C] L [ e_p ] ,  e_p \rangle 
			\sum_{q\in \N} \langle [ X_{|h|} \otimes X_{|h|} -C] L [ e_q ] ,  e_q \rangle \Big \}.
$$
To establish convergence for $n \to \infty$, we use the dominated convergence theorem. 
Dominated convergence is guaranteed by the following calculation:
\begin{align} \label{C2_1}
&\lim_{n \to \infty}  \sum_{|h|<n} \Big|  \E \bigg\{\sum_{p\in \N}\langle 
[ X_0 \otimes X_0 -C] L [ e_q ] ,  e_p \rangle \sum_{q\in \N}\langle [ X_{|h|} 
\otimes X_{|h|} -C] L [ e_q ] ,  e_q \rangle \bigg\} \Big|\\
\le & 2\sum_{h \ge 0} \Big|  \E \bigg\{\sum_{p\in \N}\big
<[ X_0 \otimes X_0 -C] L [ e_q ] ,  e_p \rangle \sum_{q\in \N}\langle 
[ X_h \otimes X_h -C] L [ e_q ] ,  e_q \rangle \bigg\} \Big|\nonumber\\
\le & \sum_{h \ge 0} 2\sqrt{\phi(h)}   \E \bigg\{ 
\Big[ \sum_{p\in \N}\langle [ X_0 \otimes X_0 -C] L [ e_q ] ,  e_p \rangle  \Big]^2\bigg\}. \nonumber
\end{align}
Here we have again used the covariance inequality from \eqref{Dehling_Lemma}. 
Given the summability of $\sqrt{\phi(h)}$ 
the first sum is finite. Turning to the expectation on the right, 
we observe that it can be expressed more compactly as 
$$
 \E  \Tr \left( \sqrt{n}(\hat{C}_n -C)L \right)^2,
$$
which is upper bounded by
\begin{equation}
 \E  \left\vvvert \sqrt{n}( \hat{C}_n - C ) L \right\vvvert_1^2. \label{applyDavies}
\end{equation}
Theorem~5.6.7 in \cite{Davies-Book-2007}
states that for a trace class operator $A$ and a bounded operator $B$ the inequality 
$\vvvert A B \vvvert_1 \le \vvvert A^* \vvvert_1 \vvvert B \vvvert_\mathcal{L}$ 
holds. Applying this inequality to  \eqref{applyDavies} yields
$$ 
\E \big\vvvert\big[\sqrt{n}(\hat{C}_n -C) \big\vvvert_1^2 \vvvert L \vvvert_\mathcal{L}^2,
$$
which is finite as can be deduced by similar reasoning as in (C1). 
Thus \eqref{C2_1} is indeed finite and the dominated convergence 
theorem ensures that $\Phi(L)$ is well defined.

In order to show condition (C4), we define
$$
F_k:= \textnormal{span}\{e_p \otimes e_q : 1 \le p, q \le k\} \subset \mathcal{S}_1(H, H),
$$
 which is a $k^2$-dimensional subspace. We then conclude that 
 \begin{align*}
&\lim_{k \to \infty} \sup_n  \E \inf_{f \in F_k}\Big\vvvert \frac{1}{\sqrt{n}}\sum_{i=1}^n \big[X_i \otimes X_i -C\big]-f\Big\vvvert_1^2 \\
\le &\lim_{k \to \infty}  \sup_n \frac{1}{n}  \E \Big\vvvert \sum_{p,q>k}e_p \otimes e_q\sum_{i=1}^n  \langle  X_i, e_p\rangle \langle  X_i, e_q\rangle -\E\langle  X_i, e_p\rangle \langle  X_i, e_q\rangle \Big\vvvert_1^2 \\
\le & 2\lim_{k \to \infty}  \sup_n \frac{1}{n}  \E \Big\vvvert \sum_{p,q>k}e_p \otimes e_q\sum_{i=1}^n  \langle  X_i, e_p\rangle \langle  X_i, e_q\rangle \Big\vvvert_1^2 \\
+&2\lim_{k \to \infty}  \sup_n \frac{1}{n}  \E \Big\vvvert \sum_{p,q>k}e_p \otimes e_q\sum_{i=1}^n   \E  \langle  X_i, e_p\rangle \langle  X_i, e_q\rangle \Big\vvvert_1^2.
\end{align*}
The first term on the right is according to Jensen's inequality upper bounding the second one. 
It will thus suffice to show that the former converges to $0$. Notice that 
$\sum_{p,q>k}e_p \otimes e_q\sum_{i=1}^n   \langle  X_i, e_p\rangle \langle  X_i, e_q\rangle $ 
is a symmetric, positive definite operator. Consequently its trace norm equals its trace. 
\begin{align*}
& \lim_{k \to \infty}  \sup_n \frac{1}{n}  \E \Big\vvvert  \sum_{p,q>k}e_p \otimes e_q\sum_{i=1}^n  \langle  X_i, e_p\rangle \langle  X_i, e_q\rangle \Big\vvvert_1^2\\
=& \lim_{k \to \infty}  \sup_n \frac{1}{n}  \E  \Tr\Big( \sum_{p,q>k}e_p \otimes e_q\sum_{i=1}^n  \langle  X_i, e_p\rangle \langle  X_i, e_q\rangle  \Big)^2\\
=&\lim_{k \to \infty}  \sup_n \frac{1}{n}  \E  \Big( \sum_{p>k}\ \sum_{i=1}^n  \langle  X_i, e_p\rangle ^2 \Big)^2 \\
=&\lim_{k \to \infty}  \sup_n \frac{1}{n}  \E   \sum_{i,j=1}^n \sum_{p>k} \langle  X_i, e_p\rangle ^2 \sum_{q>k} \langle  X_j, e_q\rangle ^2 \\
=&\lim_{k \to \infty}  \sup_n  \sum_{|h|<n} \big( 1-\frac{|h|}{n}\big) \sqrt{\phi(h)}  \E   \sum_{p>k} \langle  X_0, e_p\rangle ^2 =0
\end{align*}
In the last step we have again used the summability of the $\phi(h)$ as well as the fact that 
$$
\lim_{k \to \infty} \E  \sum_{p>k} \langle  X_0, e_p\rangle ^2=0.
$$
This concludes the proof of Theorem~\ref{Theo_CLT_EmpCovFunc}.
Theorem~\ref{Theo_funcCLT_EmpCovFunc} now follows immediately by Theorem~3.2 in \cite{Samur1987} (Condition~I). 

\hfill $\square$

%
%
%
%
%
%
%
%
%


\subsection{Proof of Theorem~\ref{Theo_FrechetDeriv} }\label{Subsec_ProofFrechDeriv} \label{secA2}

The proof of Fr\'{e}chet differentiability in $\{\lambda \mapsto \lambda C\}$ is similar for all three maps. We have to demonstrate that each of them is the composition of Fr\'{e}chet differentiable functions, which yields differentiability of the composition by the chain rule  \citep[see for instance Section~9 of Chapter~3 in][]{vdVaartWellnerBook1996}. For notational ease we subsequently do not distinguish between a function $f$ and its evaluation $f(\lambda)$. In particular we write $\lambda C$ instead of  $\{\lambda \mapsto \lambda C\}$ .\\

To show Fr\'{e}chet differentiability of $\mathbf{F}^{\rm tr}$ we notice that $\mathbf{F}^{\rm tr}$ is the 
composition of the two maps
$$
	\mathbf{A}^{(\rm tr)}: \mathcal{C} \left(  I , \mathcal{S}_1(H, H) \right) \to 
		\mathcal{C} \left( I ,\mathcal{S}_1( H_1, H_1) \times \mathcal{S}_1( H_2, H_2)  \times \mathcal{S}_1(H, H) \right) 
$$
	and 
$$
	\mathbf{B}^{(\rm tr)} : \mathcal{C} \left( I ,\mathcal{S}_1( H_1, H_1 ) \times \mathcal{S}_1( H_2, H_2)  \times \mathcal{S}_1(H, H) \right) 
		 \to 
		 \mathcal{C} \left(  I ,\mathcal{S}_1(H,H) \right) 
$$
pointwise defined as
$$
\mathbf{A}^{(\rm tr)}: L(\lambda) \mapsto \left(\begin{array}{c} 
		\operatorname{Tr}_2  [ L( \lambda ) ] / \Tr [L(\lambda ) ]
	\\ 
		  \operatorname{Tr}_1  [ L(\lambda ) ]
	\\ 
		  L(\lambda)
		  \end{array}
		  \right)
		  $$
and 
$$
	\mathbf{B}^{(\rm tr)} : \left(\begin{array}{c} 
		L_1(\lambda )
	\\ 
		 L_2(\lambda)
	\\ 
		  L_3(\lambda)
		  \end{array}
		  \right) \mapsto L_1(\lambda)\otimes L_2(\lambda)-L_3(\lambda).
$$

Beginning with $\mathbf{A}^{(\rm tr)}$, we have to show Fr\'{e}chet differentiability of each component. 
By definition  \citep[see  Section~3.9 in][]{vdVaartWellnerBook1996} the Fr\'echet derivative has to be a bounded, linear map. In the following we will use the continuity of the partial traces to find the desired Fr\'echet derivatives of the maps $\mathbf{A}^{ \rm tr }_i$, $i=1, 2, 3$.

 Due to linearity and continuity with respect to the maximum norm (defined in \eqref{maxnorm}) of the partial trace   and the identity,  the second and the third component of $\mathbf{A}^{(\rm tr)}$ are trivially
Fr\'{e}chet  differentiable, where each map is its own derivative. For the first component $\mathbf{A}^{(\rm tr)}_1$ a straightforward calculation yields that for any function $h \in \mathcal{C} \left(  I , \mathcal{S}_1(H,H) \right)$ the following equation holds pointwise, for any $\lambda \in I$ 
\begin{align*}
		\mathbf{A}_1^{(\rm tr)}[\lambda C+h(\lambda)]-\mathbf{A}_1^{(\rm tr)}[\lambda C] =\frac{
			\operatorname{Tr}_2 [C \lambda + h(\lambda) ]
		}{
			\operatorname{Tr} [ C \lambda + h(\lambda) ] 
		}
		-
		\frac{\operatorname{Tr}_2 [ C \lambda  ] }{ \operatorname{Tr} [ C \lambda ] }
	=
		\frac{ 
			\operatorname{Tr}_2 [ h(\lambda) ] \operatorname{Tr} [ C] 
			- 
			\operatorname{Tr}_2[ C ] \operatorname{Tr} [ h(\lambda) ]
		}{
			\operatorname{Tr} [ C\lambda + h ( \lambda ) ] \operatorname{Tr} [ C ]
		}.
\end{align*}
We hence define the derivative $D\mathbf{A}_1^{(\rm tr)}(\lambda C)$ of $\mathbf{A}_1^{(\rm tr)}$ in $\lambda C$ as follows
$$ 
	D \mathbf{A}_1^{(\rm tr)} ( \lambda C ) : 
		h(\lambda)
		 \mapsto 
		\frac{ 
			\operatorname{Tr}_2 [ h(\lambda) ] \operatorname{Tr} [C]
			- 
			\operatorname{Tr}_2[ C ] \operatorname{Tr} [ h(\lambda) ]
		}{
		\lambda \operatorname{Tr} [C ]^2
		}.
$$
This function is linear and bounded, since $\lambda$ is 
bounded away from $0$. We now have to prove that this function is indeed 
the Fr\'{e}chet derivative. Let $h \in  \C ( I, \mathcal{S}_1(H,H) )$ be a function. Some simple calculations show that
\begin{align*}
&
	\mathbf{A}_1^{(\rm tr)}[\lambda C+h(\lambda)] 
	-
	 \mathbf{A}_1^{(\rm tr)}[\lambda C]
	 -
	 D \mathbf{A}_1^{(\rm tr)}  (   \lambda C ) [h(\lambda)]
	\\[1ex]
=&
	\frac{ 
			\operatorname{Tr}_2 [ h(\lambda) ]\operatorname{Tr} [ C ] 
			- 
			\operatorname{Tr}_2 [ C ] \operatorname{Tr} [ h(\lambda) ]
		}{
			\operatorname{Tr} [ C\lambda + h ( \lambda ) ] \operatorname{Tr} [ C ]
		}
	-
	\frac{ 
			\operatorname{Tr}_2 [ h(\lambda) ] \operatorname{Tr} [ C ]
			- 
			\operatorname{Tr}_2 [C ] \operatorname{Tr} [ h(\lambda) ]
		}{
		\lambda \operatorname{Tr} [ C ]^2
		}
	\\[1ex]
=& 
	\frac{ - \operatorname{Tr}_2 [ h(\lambda) ] \operatorname{Tr} [ h ( \lambda ) ] 
		}{
		\lambda \operatorname{Tr} [ C\lambda + h ( \lambda ) ] \operatorname{Tr} [ C ] }
	+
	\frac{  \operatorname{Tr} [ h(\lambda) ]^2 \operatorname{Tr}_2 [ C ]  
		}{
		\lambda \operatorname{Tr} [ C\lambda + h ( \lambda ) ] \operatorname{Tr} [ C ]^2 }
\end{align*}
Noting that for $ \lambda \in I$, both $1 / \lambda$ and $1/ \vvvert \operatorname{Tr} [ C\lambda] \vvvert_1$ are bounded uniformly in $\lambda$ by a constant and recalling that the trace and the partial traces are Lipschitz continuous with respect to the trace norm, it  follows that the above expression is of the order $\mathcal{O}\big( \|h\|_\infty^2\big)$.
Hence $D  \mathbf{A}_1^{(\rm tr)} ( \lambda C ) $ indeed is the Fr\'{e}chet derivative of $\mathbf{A}_1^{(\rm tr)}$ in $ \lambda C$.  
Summarizing we have that
	$$
		D \mathbf{A}^{(\rm tr)} (  \lambda C )  : 
		L( \lambda ) 
	\mapsto
		  		\left(\begin{array}{c} 
			
			\Tr_2 [ L (\lambda) ] /  \Tr [ \lambda C ]
			- 
			\Tr [ L( \lambda)] \Tr_2[ C] / (  \lambda \Tr [ C ]^2 )
		\\ 
			 \Tr_1[ L(\lambda) ]
		\\ 
			    L(\lambda)
		  \end{array}  \right) .
$$

The map $\mathbf{B}^{(\rm tr)}  $ is differentiable in 
$\mathbf{A}^{(\rm tr)} [C \lambda] = ( C_1 ( \lambda ), C_2 ( \lambda ) ,  \lambda C )$, with derivative, given as 
	$$
		D \mathbf{B}^{(\rm tr)} (  C_1 ( \lambda ), C_2 ( \lambda ) ,  \lambda C ) : 
		\left(\begin{array}{c} 
			L_1(\lambda )
		\\ 
			 L_2(\lambda)
		\\ 
			  L_3(\lambda)
		  \end{array}  \right) 
	\mapsto
		  C_1 ( \lambda ) \otimes L_2( \lambda )
		  +
		  L_1(\lambda)\otimes C_2 ( \lambda )  
		  -
		  L_3(\lambda).
$$
The proof is similar, but easier than the one for the differential of $\mathbf{B}^{(\rm prod)}$ and hence omitted.
By the above arguments $\mathbf{F}^{(\rm tr)}$ is differentiable in  $\lambda C$.\\

We now turn to the differentiability of the map $\mathbf{F}^{(\rm prod)}$. We can again decompose this map into two simpler ones

$$
	\mathbf{A}^{(\rm prod)}: \mathcal{C} \left(  I , \mathcal{S}_2(H, H) \right) \to 
		\mathcal{C} \left( I ,\mathcal{S}_2( H_1, H_1 ) \times \mathcal{S}_2( H_2, H_2)  \times \mathcal{S}_2(H,H) \right) 
$$
	and 
$$
	\mathbf{B}^{(\rm prod)} : \mathcal{C} \left(I,\mathcal{S}_2( H_1, H_1) \times \mathcal{S}_2( H_2, H_2)  \times \mathcal{S}_2(H, H) \right) 
		 \to 
		 \mathcal{C} \left(  I ,\mathcal{S}_2(H, H) \right) 
$$
pointwise defined as
$$
\mathbf{A}^{(\rm prod)}: L(\lambda) \mapsto \left(\begin{array}{c} 
		P_2(L(\lambda), \Delta_2)
	\\ 
		 P_1\big(L(\lambda), P_2(L(\lambda), \Delta_2)\big)
	\\ 
		  L(\lambda)
		  \end{array}
		  \right)
		  $$
and 
$$
	\mathbf{B}^{(\rm prod)} : \left(\begin{array}{c} 
		L_1(\lambda )
	\\ 
		 L_2(\lambda)
	\\ 
		  L_3(\lambda)
		  \end{array}
		  \right) \mapsto \frac{L_1(\lambda)\otimes L_2(\lambda)}{\vvvert L_1(\lambda) \vvvert_2^2}-L_3(\lambda).
$$
We begin with the differentiability of $\mathbf{A}^{(\rm prod)}$. As in the previous case one sees that the first and the last component are obviously differentiable due to linearity and boundedness. For the component map $\mathbf{A}_2^{(\rm prod)}$, by bilinearity of the maps $P_1, P_2$, the difference $\mathbf{A}_2^{(\rm prod)}[\lambda C+h(\lambda)]-\mathbf{A}_2^{(\rm prod)}[\lambda C]$ equals
$$ 
		P_1( \lambda C, P_2(h(\lambda), \Delta_2) )
	+ 
		P_1( h(\lambda), P_2(\lambda C, \Delta_2) )
	+
		P_1( h(\lambda), P_2( h(\lambda), \Delta_2) ).
$$
Due to the bilinearity of the partial products, the first and the second term are linear in $h(\lambda)$. The third term is obviously of order $\mathcal{O}(\|h\|_\infty^2)$. 
Hence the derivative of $\mathbf{A}_2^{(\rm prod)}$ in $\lambda C$ is pointwise defined as 
$$
		D \mathbf{A}_2^{(\rm prod)} ( \lambda C ) [h(\lambda)]
	:= 
		P_1 ( \lambda C, P_2( h(\lambda), \Delta_2 ) ) + P_1( h(\lambda), P_2( \lambda C, \Delta_2) ).
$$
We claim that the derivative of $\mathbf{B}^{(\rm prod)}$ in any function $T(\lambda):=(T_1 (\lambda ), T_2 (\lambda ), T_3 (\lambda ))$ 
is given by the bounded linear map
	$$
		D \mathbf{B}^{(\rm prod)}(T(\lambda))  : 
		\left(\begin{array}{c} 
			L_1(\lambda )
		\\ 
			 L_2(\lambda)
		\\ 
			  L_3(\lambda)
		  \end{array}  \right) 
	\mapsto
		  \frac { 
		  	L_1(\lambda ) \otimes T_2(\lambda)  + T_1(\lambda ) \otimes L_2(\lambda) 
			}{ \vvvert T_1 ( \lambda ) \vvvert_2^2 }
		   +
		    \frac { 
		    	- 2 \langle L_1(\lambda ) , T_1(\lambda ) \rangle ( T_1(\lambda ) \otimes T_2 (\lambda ) )
		    	}{ \vvvert T_1 ( \lambda ) \vvvert_2^4 }
		    -
		    L_3( \lambda ).
$$
We sketch the arguments of the proof.
First note that 
\begin{align*}
&	\mathbf{B}^{(\rm prod)}
		[ ( L_1 (\lambda ), L_2 (\lambda ), L_3 (\lambda ) ) + ( h_1 ( \lambda ), h_2 ( \lambda ), h_3 ( \lambda ) ) ]
	-
	\mathbf{B}^{(\rm prod)}[ ( L_1 (\lambda ), L_2 (\lambda ), L_3 (\lambda ) ) ]
\\[1ex]
=&
	- h_3 ( \lambda ) 
	+
	\frac{
	h_1 ( \lambda ) \otimes L_2 (\lambda ) + L_1 ( \lambda ) \otimes h_2 (\lambda ) 
	}{
	\vvvert L_1 ( \lambda ) + h_1 ( \lambda ) \vvvert_2^2
	}
	-
	\frac{
	2 \langle h_1 ( \lambda ), L_1 ( \lambda ) \rangle L_1 ( \lambda )  \otimes L_2 ( \lambda ) 
	}{
	\vvvert L_1 ( \lambda ) + h_1 ( \lambda ) \vvvert_2^2 \vvvert L_1 ( \lambda )  \vvvert_2^2
	}
	+
	\mathcal{O}\big( \|h\|_\infty^2\big).
\end{align*}
Subtracting $D \mathbf{B}^{(\rm prod)}  ( T(\lambda)) 
[ ( h_1 ( \lambda ), h_2 ( \lambda ), h_3 ( \lambda ) ) ]$ from the above yields a remainder, which is of order $\mathcal{O}\big( \|h\|_\infty^2\big)$, as a short calculation reveals. Hence $D \mathbf{B}^{(\rm prod)}  ( T(\lambda)) $ must be the Fr\'{e}chet differential of $\mathbf{B}^{(\rm prod)}$ in $T(\lambda) $.\\

Finally we turn to the differentiability of $\mathbf{F}^{(\rm opt)}$. Let $\{u_i\}_{ i \ge 1 }$ be a basis of the space $H_1$ and $\{v_j\}_{ j \ge 1 }$ a basis of $H_2$. Then $\{u_i \otimes v_j\}_{i,j \ge 1 }$ is a tensor basis of $H$ and $\{u_i \otimes v_j \otimes u_k \otimes v_l\}_{i,j,k,l \ge 1}$ of $H \otimes H$. Recall that $H \otimes H \cong \mathcal{S}_2(H, H)$, i.e. we may identify any operator $L \in \mathcal{S}_2(H, H)$, with the tensor basis expansion
$$ L = \sum_{i,j,k,l} \alpha_{i,j,k,l} u_i \otimes v_j \otimes u_k \otimes v_l,$$
where $ \alpha_{i,j,k,l}  \in \mathbb{R}$. Similarly any sequential operator $L(\lambda)$  in the space $\mathcal{C}(I, \mathcal{S}_2(H, H))$ can be expressed  by a tensor expansion, where accordingly the coefficients $\alpha_{i,j,k,l}(\lambda)$ depend on the sequential parameter. We can hence define the map
$$\mathbf{A}^{(\rm opt)}: \begin{cases} 
     \mathcal{C}(I, \mathcal{S}_2(H, H)) \to  \mathcal{C}(I, \mathcal{S}_2(H_1 \otimes H_1, H_2 \otimes H_2)) \\ 
       \sum_{i,j,k,l} \alpha_{i,j,k,l}(\lambda) u_i \otimes v_j \otimes u_k \otimes v_l \mapsto \sum_{i,j,k,l} \alpha_{i,j,k,l}(\lambda) u_i \otimes  u_k  \otimes v_j \otimes v_l  \\ 
   \end{cases}  $$
which is bijective, linear and norm preserving (in particular bounded).  Notice that $\mathbf{A}^{(\rm opt)}$ is a version of the map $\Pi$, introduced in Section \ref{Subsubsec_OptAppr} for sequential operators. Secondly, denote by $L(\lambda)^*$ the adjoint and consider the map
 $$
 	\mathbf{B}^{(\rm opt)}: \begin{cases} 
		\mathcal{C}(I, \mathcal{S}_2(H_1 \otimes H_1, H_2 \otimes H_2)) 
		\to 
		\mathcal{C}(I, \mathcal{S}_2^{\rm sym}(H_1 \otimes H_1, H_1 \otimes H_1) 
				\times 
				\mathcal{S}_2^{\rm sym}(H_2 \otimes H_2, H_2 \otimes H_2) )
\\ 
	 	 L(\lambda) 
	 \mapsto 
	 	\left( L( \lambda )^* L( \lambda ), L( \lambda ) L( \lambda )^* \right)\\ 
   \end{cases}.
$$
  Here the space $\mathcal{S}_2^{\rm sym}(H_i \otimes H_i, H_i \otimes H_i)$  is defined as the vector space of all Hilbert--Schmidt operators in $\mathcal{S}_2(H_i \otimes H_i, H_i \otimes H_i)$, which are symmetric. The inner product on this space is  induced by that on $\mathcal{S}_2(H_i \otimes H_i, H_i \otimes H_i)$. The map $\mathbf{B}^{(\rm opt)}$ is Fr\'{e}chet differentiable in the sequential operator $A_0(\lambda):=\mathbf{A}^{(\rm opt)}[\lambda C]$ and its derivative 
  $$
  	D \mathbf{B}^{(\rm opt)} ( A_0(\lambda) )
	: L(\lambda) 
	\mapsto 
	\left(	
		L^*(\lambda) L_0(\lambda) + L_0^*(\lambda) L(\lambda) 
		,
		L (\lambda)L_0^*(\lambda) + L_0(\lambda) L^*(\lambda)
	\right).
$$ 
Note that the image of $D \mathbf{B}^{(\rm opt)} ( A_0(\lambda))$ consists of tupels of pointwise symmetric operators. Recall that for symmetric operators both eigenvectors and eigenvalues are well defined objects. 
We hence define for $i=1,2$ the maps
   $$
       \mathbf{C}_i^{(\rm opt)}:
    \begin{cases}
            \mathcal{C}(I, \mathcal{S}_2^{\rm sym}( H_i \otimes H_i, H_i \otimes H_i ) )
         \to
         \mathcal{C}(I,H_i \otimes H_i )\\
            L(\lambda) \mapsto v_1^{ L ( \lambda ) } ,\\
   \end{cases}
$$
   where $v_1^{ L ( \lambda ) } $ denotes the normalized eigenvector of $L(\lambda)$ corresponding to the largest eigenvalue. Let $B_0(\lambda):= ( \mathbf{B}^{\rm (opt)} \circ \mathbf{A}^{\rm (opt)} ) [\lambda C] $, which has unequal first two eigenvalues for all $\lambda \in I$.  This assertion holds because of Assumption~\ref{assumption_1} (4), which implies $\gamma_1>\gamma_2$ for the first two singular values of $\Pi[C]=\mathbf{A}^{\rm (opt)}[C]$.  $\mathbf{C}_i$ is differentiable in $B_0 (\lambda)$ with derivative
  $$
  	D \mathbf{C}_{i}^{(\rm opt)} ( B_0 (\lambda) ): 
	h(\lambda) 
	\mapsto 
	\sum_{k >1} \frac{\langle h(\lambda),  v_1^{ B_0(\lambda)} \otimes v_k^{ B_0 (\lambda)} \rangle v_k^{ B_0 (\lambda)}  }{\gamma_1^{ B_0(\lambda) } - \gamma_k^{ B_0 (\lambda) }  } ,
$$
where $\gamma^{L ( \lambda ) }_k$ denotes the $k$-th largest eigenvalue of the operator $L(\lambda)$ and $v_k^{ L ( \lambda ) }$ the corresponding eigenvector. Analogously we observe that the map
$$ \mathbf{C}_3^{(\rm opt)}: \begin{cases} 
    \mathcal{C}(I, \mathcal{S}_2^{\rm sym}(H_1 \otimes H_1, H_1 \otimes H_1) ) \to \mathcal{C}(I,\mathbb{R})\\ 
     L(\lambda) \mapsto \gamma_1^{L(\lambda)},\\ 
   \end{cases} $$
is differentiable in $B_0(\lambda)$ with derivative
$$ 
	D \mathbf{C}_{3}^{(\rm opt)} ( B_0(\lambda) )
	: h(\lambda) 
	\mapsto 
	\langle h(\lambda), v_1^{B_0(\lambda)} \otimes v_1^{B_0(\lambda)} \rangle. 
$$
The proof of the differentiability of these maps is because of its techincal nature  deferred to the end of this proof.
By differentiability of each component, the combined map  
 $$ 
 	\mathbf{C}^{(\rm opt)}: 
	\begin{cases} 
   			 \mathcal{C}(I, \mathcal{S}_2^{\rm sym}(H_1 \otimes H_1, H_1 \otimes H_1) 
			\times 
			\mathcal{S}_2^{\rm sym}(H_2 \otimes H_2, H_2 \otimes H_2) )
		 \to 
		 	\mathcal{C}(I,   ( H_1 \otimes H_1 ) \times ( H_2 \otimes H_2 ) \times \mathbb{R} )
		 
	\\ 
    		\left( L_1( \lambda ), L_2 ( \lambda ) \right)  
		 \mapsto 
		  \left( 
		  	\mathbf{C}^{(\rm opt)}_1 ( L_1( \lambda ) ), 
			\mathbf{C}^{(\rm opt)}_2 ( L_2 ( \lambda ) ),
			\mathbf{C}^{(\rm opt)}_3  ( L_1( \lambda ) )
		  \right)
	\\ 
   \end{cases} 
$$
is differentiable.
 Next, the map 
 $$ 
 	\mathbf{D}^{(\rm opt)}: 
	\begin{cases} 
   		 \mathcal{C}(I,  ( H_1 \otimes H_1 ) \times ( H_2 \otimes H_2 )  \times \mathbb{R} )
		 \to 
		  \mathcal{C}(I, \mathcal{S}_2(H_1 \otimes H_1, H_2 \otimes H_2)) 
	\\ 
    		\left( V( \lambda ), W( \lambda ) ,  r( \lambda ) \right)  
		 \mapsto 
		 r(\lambda) \lambda [V(\lambda) \otimes W(\lambda)]
	\\ 
   \end{cases} 
  $$
   is Fr\'{e}chet differentiable in any $3-$tuple $( V_0(\lambda), W_0(\lambda), r_0(\lambda) )$ with derivative 
   $D\mathbf{D}^{(\rm opt)} ( V_0(\lambda), W_0(\lambda), r_0(\lambda) ) $ given as
  $$ 
		 \left(  V( \lambda ), W( \lambda ), r( \lambda ) \right)
		 \mapsto 
		 \lambda \left\{ r_0(\lambda)  [V_0(\lambda) \otimes W(\lambda)]+r(\lambda) [V_0(\lambda) \otimes W_0(\lambda)]+r_0(\lambda) [V(\lambda) \otimes W_0(\lambda)] \right\},
$$
   which follows by elementary calculations.
Finally $(\mathbf{A}^{(\rm opt)})^{-1}$  is bounded, linear and thus differentiable. 
We now conclude that the map
$$\mathbf{E}^{(\rm opt)}:=   (\mathbf{A}^{(\rm opt)})^{-1} \circ \mathbf{D}^{(\rm opt)} \circ \mathbf{C}^{(\rm opt)} \circ \mathbf{B}^{(\rm opt)} \circ \mathbf{A}^{(\rm opt)}:   \mathcal{C}(I, \mathcal{S}_2(H, H)) \to \mathcal{C}(I, \mathcal{S}_2(H, H)),$$
which maps a sequential operator $L(\lambda)$ to its optimal separable approximation $L^{\rm opt}(\lambda)$
is differentiable in $\lambda \mapsto \lambda C$. Hence it follows immediately that $\mathbf{F}^{(\rm opt)}$ also is.\\

We conclude this section by a proof of the differentiability of the component maps $\mathbf{C}^{(\rm opt)}_1, \mathbf{C}^{(\rm opt)}_2$ and $\mathbf{C}^{(\rm opt)}_3$.
We first make an observation, that will be used below.
For an operator
$
A 
 \in
 \mathcal{ C }( I,  \mathcal{S}_2^{\rm sym}(H_i \otimes H_i, H_i \otimes H_i) )$, $i=1,2$, by the spectral theorem we have the representation
$$ A  ( \lambda )
=
 \sum_{ k \ge 1 } \gamma_k^{ A ( \lambda ) }  v_k^{ A ( \lambda ) }  \otimes v_k^{ A ( \lambda ) }, 
 $$where $\gamma_k^{ A ( \lambda ) }  \in C( I, \mathbb{ R } )$, $v_k^{ A ( \lambda ) }  \in \mathcal{ C } ( I, H_i \otimes H_i )$, $i=1,2$, are the eigenvalues and  respective normalized eigenvectors of
  $A ( \lambda ) $.
Combing this identity with
$
 \sum_{ j \ge 1 } \langle  v_j^{ A ( \lambda ) }   , x  \rangle  v_j^{ A ( \lambda ) }
  =
  x$, for $ x \in  H_i \otimes H_i $ we see that
$$
        (  \gamma_1^{ A ( \lambda ) }  \operatorname{Id}  - A ( \lambda ) ) [ x  ]
    =
         \sum_{ i \ge 2 } ( \gamma_1^{ A ( \lambda ) }  - \gamma_i^{ A ( \lambda ) }  )
          \langle  v_i^{ A ( \lambda ) }   , x  \rangle  v_i^{ A ( \lambda ) }.
$$
We further define the symmetric operator
$
B
\in
 \mathcal{ C }( I,  \mathcal{S}_2^{\rm sym}(H_i \otimes H_i, H_i \otimes H_i) )
 $ as
$$
    B ( \lambda )
        :y
    \mapsto
    \sum_{ j \ge 2 }
        \frac{ \langle v_j^{ A ( \lambda ) } , y \rangle
        }{
           \gamma_1^{ A ( \lambda ) }  -  \gamma_j^{ A ( \lambda ) }   }
        v_j^{ A ( \lambda ) }.
$$
Now a simple calculation shows that
$B ( \lambda )
    [ ( \gamma_1^{ A ( \lambda ) }  \operatorname{Id}  - A ( \lambda ) ) [ x ] ]
= x$,
for $x \in \operatorname{Ker} ( \gamma_1^{ A ( \lambda ) }  \operatorname{Id}  - A ( \lambda ) )^{ \bot }$. Hence $B ( \lambda )$ is the left inverse of
$( \gamma_1^{ A ( \lambda ) }  \operatorname{Id}  - A ( \lambda )  )$ on this space and by abuse of notation we will denote it
by $ ( \gamma_1^{ A ( \lambda ) }  \operatorname{Id}  - A ( \lambda ) )^{ - 1 }$.

We first consider the eigenvectors and only give the proof for $\mathbf{C}_1^{(\rm opt)}$, as the case of
$\mathbf{C}_2^{(\rm opt)}$ is analogous. Let
$A, h \in \mathcal{C}(I, \mathcal{S}_2^{\rm sym}(H_1 \otimes H_1, H_1 \otimes H_1) ) $ and $\|h \|_\infty\to 0 $, then
\begin{align*}
        v_1^{ ( A  + h ) ( \lambda ) }   - v_1^{ A ( \lambda)   }
    =&
        \frac{ 1 }{  \gamma_1^{ ( A  + h ) ( \lambda ) }  }
        \left(
         ( A + h ) ( \lambda ) [ v_1^{ ( A  + h ) ( \lambda ) }  ] -  \gamma_1^{ ( A  + h ) ( \lambda ) }  v_1^{ A  ( \lambda)  }
        \right)
    \\
    =&
        \frac{ 1 }{   \gamma_1^{ ( A  + h ) ( \lambda ) } }
        \left(
        A ( \lambda ) [ v_1^{ ( A  + h ) ( \lambda ) }  - v_1^{ A  ( \lambda ) }   ]
        +
         h ( \lambda ) [ v_1^{ ( A  + h ) ( \lambda ) }   ]
        +
         v_1^{ A ( \lambda )   }
        \left(
             \gamma_1^{ A ( \lambda )  }   -  \gamma_1^{ ( A  + h ) ( \lambda ) }
        \right)
        \right)
\end{align*}
where we used the properties of eigenvalues and eigenvectors.
We further simplify by multiplying both sides by $\gamma_1^{ ( A  + h ) ( \lambda ) } $ and rearranging terms, which gives
$$
    \left(  \gamma_1^{ ( A  + h ) ( \lambda ) }  \operatorname{Id} - A ( \lambda ) \right)
    \left[ v_1^{ ( A  + h ) ( \lambda ) } - v_1^{ A  ( \lambda ) }  \right]
=
    h( \lambda ) [ v_1^{ ( A  + h ) ( \lambda ) } ]
    +
    v_1^{ A( \lambda)  }
     \left(
             \gamma_1^{ A ( \lambda )  }   -  \gamma_1^{ ( A  + h ) ( \lambda ) }
    \right)
$$
Define the function
$$
    D \mathbf{C}_1^{(\rm opt)}  ( A  ( \lambda ) ) [ h ( \lambda ) ]
:=
    ( \gamma_1^{ A ( \lambda )  } \operatorname{Id} - A ( \lambda ) )^{ - 1 }
    \left[  h ( \lambda ) [  v_1^{ A( \lambda )   }  ] \right]
$$

A calculation shows that,
\begin{align*}
& \left( v_1^{ ( A  + h ) ( \lambda ) }   - v_1^{ A ( \lambda ) }  \right) / \| h \|_\infty \\
= &    \left(  \gamma_1^{ ( A + h) ( \lambda ) }  \operatorname{Id} - A ( \lambda ) \right)^{ - 1 }
    \Big [
        \frac{ h ( \lambda ) }{ \| h \|_\infty } \left[ v_1^{ ( A  + h ) ( \lambda ) }  \right]
        +
        \frac{ v_1^{ A ( \lambda )   } }{ \| h \|_\infty }
             \left( \gamma_1^{ A ( \lambda )   }   -  \gamma_1^{ ( A  + h ) ( \lambda ) }  \right)
    \Big ]
\end{align*}

Note from the above observation that $(  \gamma_1^{ A ( \lambda ) }  \operatorname{Id}  - A ( \lambda ) )^{ - 1 } $ is bounded and
$(   \gamma_1^{ A ( \lambda ) }  \operatorname{Id}  - A ( \lambda ) )^{ - 1 } [ v_1^{ A ( \lambda ) } ] =  0$.
Using these facts and subtracting $D \mathbf{C}_1^{(\rm opt)}  ( A ( \lambda )  ) [ h ( \lambda ) ]  / \| h \|_\infty$ from the above expression yields a term, which has norm  of order
\begin{align*}
&
    \mathcal{O}
    \Big(
        \left\|
                \left(  \gamma_1^{ A + h }  \operatorname{Id} - ( A + h ) \right)^{ - 1 }
         \right\|_\infty
          \|  v_1^{ A    + h}- v_1^{ A }  \|
    \\
+&
      2 \left\|
                \left(  \gamma_1^{ A + h }  \operatorname{Id} - ( A + h ) \right)^{ - 1 }
                -
                \left(  \gamma_1^{ A   }  \operatorname{Id} - A \right)^{ - 1 }
            \right\|_\infty
     \Big),
\end{align*}
where we used the inequality $ | \gamma_1^{ A   } ( \lambda )  -  \gamma_1^{ A  + h } ( \lambda ) |
\le \vvvert ( A + h ) ( \lambda ) - A ( \lambda ) \vvvert_2$ for every $\lambda$.
Now both terms converge to 0, by continuity of the eigenvalues and eigenvectors (see Lemmas~2.2. and 2.3 in \cite{HorKokBook2012}).

Next, we derive the Fr\'{e}chet derivative of the eigenvalue map. Let again $h$ be a symmetric operator, converging to 0. Then, we have the following:
\begin{align*}
        \gamma_1^{ ( A  + h ) ( \lambda ) }  - \gamma_1^{ A ( \lambda ) }
    =&
        \langle ( A  + h ) ( \lambda ) [ v_1^{ ( A  + h ) ( \lambda )  }  ],  v_1^{ ( A  + h ) ( \lambda ) }   \rangle
        -
        \langle  A ( \lambda ) [ v_1^{ A ( \lambda ) } ], v_1^{ A ( \lambda ) } ( \lambda) \rangle
    \\
    =&
        \langle
            A  ( \lambda ) [ v_1^{ ( A  + h ) ( \lambda ) } -  v_1^{ A ( \lambda ) }  ],
             v_1^{ ( A  + h ) ( \lambda ) } +  v_1^{ A ( \lambda ) }
         \rangle
    \\
    &
    +
        \langle  h ( \lambda ) [ v_1^{ ( A  + h ) ( \lambda ) }  ]  , v_1^{ ( A  + h ) ( \lambda ) }   \rangle
\end{align*}
where $v_1^{ L  } ( \lambda) $ denotes the normalized eigenvector belonging to the largest eigenvalue of
the symmetric operator $ L ( \lambda )$.
Using the symmetry of $ A ( \lambda )$ and subtracting $\langle  h ( \lambda ) [ v_1^{ A ( \lambda ) } ], v_1^{ A ( \lambda ) }  \rangle$ we see the above expression equals
\begin{equation} \label{C_1_difference}
    \langle
              v_1^{ ( A  + h ) ( \lambda ) } -  v_1^{ A ( \lambda ) },
             A  ( \lambda ) [ v_1^{ ( A  + h ) ( \lambda ) } ]
         \rangle
         +
         \gamma_1^{ A ( \lambda ) }
         \langle
              v_1^{ ( A  + h ) ( \lambda ) } -  v_1^{ A ( \lambda ) },
               v_1^{ A ( \lambda ) }
         \rangle.
\end{equation}
Notice that according to the parallelogram law for two vectors $v, w$ of unit length the following identity holds:
$$\langle v-w, v\rangle = -\frac{1}{2} \|v-w\|^2.$$
As a consequence the second term in \eqref{C_1_difference} equals $-\gamma_1^{ A ( \lambda ) }\|v_1^{ ( A  + h ) ( \lambda ) } -  v_1^{ A ( \lambda ) }\|^2/2 =\mathcal{O}(\|h\|_\infty^2)$. 
For the first term a simple calculation yields
\begin{align*}
  & \langle
              v_1^{ ( A  + h ) ( \lambda ) } -  v_1^{ A ( \lambda ) },
             A  ( \lambda ) [ v_1^{ ( A  + h ) ( \lambda ) } ]
         \rangle \\
         =&    \langle
              v_1^{ ( A  + h ) ( \lambda ) } -  v_1^{ A ( \lambda ) },
             A  ( \lambda ) [ v_1^{ ( A  + h ) ( \lambda ) }- v_1^{ ( A   ) ( \lambda ) }]\rangle +  \gamma_1^{A(\lambda)} \langle
              v_1^{ ( A  + h ) ( \lambda ) } -  v_1^{ A ( \lambda ) }, v_1^{ ( A   ) ( \lambda )}\rangle.
\end{align*}
Both terms on the right are of order $\mathcal{O}(\|h\|_\infty^2)$, which can be seen for the first one,  by an application of Cauchy-Schwarz and for the second by using the parallelogram law.

%
%
%
%
%
%
%
%
%

\subsection{Proof of Corollary \ref{cor_main}}

Let $\rm x \in \{tr, prod, opt\}$.  We define 
$$
	H_{n}^{\rm x}(\lambda)=  M^{\rm x}[ \hat C_n(\lambda)]- M^{\rm x}[\lambda C] .
$$
In order to prove the weak convergence
\begin{equation} \label{conv_to_BM}
\{ \sqrt{n} H_{n}(\lambda) \}_{\lambda \in I}  \stackrel { d } { \to } \{ \sigma^{\rm x} \lambda\mathbb{B}(\lambda)\}_{\lambda \in I},
\end{equation}
we rewrite $H_{n}^{\rm x}(\lambda)$ as follows:
\begin{align*}
	H_{n}(\lambda) =  &	 \vvvert    \hat{C}_n(\lambda) -   \hat{C}_n^{\rm x}(\lambda) \vvvert_2^2
		-
		\lambda^{2} \vvvert    C - C^{\rm x} \vvvert_2^2
		\\
 	=&
		\Big \langle \hat{C}_n(\lambda) -  \hat{C}_n^{\rm x}(\lambda)- \lambda C + \lambda C^{\rm x} ,
 \hat{C}_n(\lambda) - \hat{C}_n^{\rm x}(\lambda) +\lambda C - \lambda C^{\rm x}  \Big  \rangle_2.
 \\
 = & 		\Big \langle \mathbf{F}^{\rm x}[\hat{C}_n(\lambda)] - \mathbf{F}^{\rm x}[C \lambda] ,
 \mathbf{F}^{\rm x}[\hat{C}_n(\lambda)]  +\mathbf{F}^{\rm x}[C \lambda]  \Big  \rangle_2.
\end{align*}
We recall Theorems~\ref{Theo_funcCLT_EmpCovFunc} and \ref{Theo_FrechetDeriv}. These imply  that
$$ \mathbf{F}^{\rm x}[\hat{C}_n(\lambda)]  +\mathbf{F}^{\rm x}[C \lambda] \stackrel{\mathbb{P}}{\to}  2 \lambda\mathbf{F}^{\rm x}[C ]$$
and by application of the functional-$\Delta$-method
$$\sqrt{n} \Big( \mathbf{F}^{\rm x}[\hat{C}_n(\lambda)]  -\mathbf{F}^{\rm x}[C \lambda]  \Big) \stackrel{d}{\to}  \mathbb{B}_G(\lambda)  ,$$
where $\mathbb{B}_G$ is a Brownian motion corresponding to the Gaussian process $G$, from Theorem \ref{Theo_funcCLT_EmpCovFunc}. We hence observe that the weak convergence
\begin{align*}
\{\sqrt{n} H_{n}(\lambda) \}_{\lambda \in I}\stackrel{d}{\to} \{2 \lambda \langle \mathbf{F}^{\rm x}[C ], \mathbb{B}_G(\lambda) \rangle_2\}_{\lambda \in I} \stackrel{d}{=} \{\sigma^{\rm x} \lambda \mathbb{B}(\lambda)\}_{\lambda \in I}, 
\end{align*}
holds, where 
\begin{equation} \label{sigmax}
(\sigma^{\rm x})^2 := 4 \mathbb{E}\big[\langle \mathbf{F}^{\rm x}[C ], \mathbb{B}_G(\lambda) \rangle_2^2\big]
\end{equation}

\hfill $\square$

\subsection{ Proof of Theorem~\ref{thmmain} and Theorem \ref{thmmain_rel}}

Both results are proved by the same arguments and for the sake of brevity we restrict ourselves to the proof of  Proof of Theorem~\ref{thmmain} .
 According to Corollary~\ref{cor_main} for $\rm x \in \{tr, prod, opt\}$ it holds that 
$$
		\left\{ \sqrt{N}\big( M^{\rm x}[ \hat C_n(\lambda)]- M^{\rm x}[\lambda C] \right) \}_{\lambda \in I }
	\stackrel{d}{\to}  
		\left\{\sigma^{\rm x}  \lambda \mathbb{B}( \lambda ) \right\}_{\lambda \in I }. 
$$
Recalling the definition of $\hat V_n^{\rm x}$ in \eqref{Eq_EstimVariance-L2norm} we observe that
\begin{align*}
			&\Big[ \int_I \big(M^{\rm x}[\hat C_n( \lambda )] - \lambda^2M^{\rm x}[\hat C_n)]  \big )^2 d\nu(\lambda)\Big]^{1/2}\\
			= &\Big[ \int_I \big(\big\{ M^{\rm x}[\hat C_n( \lambda )] -M^{\rm x}[ C \lambda ] \big\} + \big\{ \lambda^2 M^{\rm x}[ C  ]  -\lambda^2M^{\rm x}[\hat C_n)] \big\} \big )^2 d\nu(\lambda) \Big]^{1/2}\\
			\stackrel{d}{\to} & \sigma^{\rm x} \Big[  \int_I \lambda^2 \big(\mathbb{B}(\lambda)-\lambda \mathbb{B}(1) \big)^2 d\nu(\lambda) \Big]^{1/2}.
\end{align*}
By the continuous mapping theorem weak convergence of the joint vector follows 
$$
	 	\sqrt{N} \left( M^{\rm x}[ \hat C_n] - M^{\rm x}[C], \hat V_n^{\rm x} \right) 
	 \stackrel{d}{\to} 
	 	\sigma^{\rm x} 
		\Big ( 
			\mathbb{B}(1)
			, 
			\Big [\int_I \lambda^2 \left( \mathbb{B}(\lambda) -\lambda \mathbb{B}(1) \right)^2 d\nu(\lambda) \Big ]^{1/2} \Big ) ,
$$
which entails the weak convergence in \eqref{hd3} of the self-normalized statistic and proves the result. \hfill $\square$

%
%
%
%
%
%
%
%
%

%
%
%
%
%
%
%
%
%

\subsection{ Proof of Theorem \ref{thm1}}

We will restrict ourselves to a proof of the result  in \eqref{hd3} for the quantity $ M^{\rm x}$. The corresponding statement for the relative deviation   $ M_{{\rm rel}}^{\rm x}$ in 
\eqref{hd4}  follows by the same arguments.  
Recalling the definition of $\hat I_n^{\rm x}$ in \eqref{hd5} we obtain by an elementary calculation that 
	$$
		\mathbb{P}  \{ M^{\rm x}[C] \in \hat I_n^{\rm x} \}
	= 
		\mathbb{P} \{ q_{\alpha/2} \hat V_n^{\rm x} \le \hat M^{\rm x}[\hat C_n] 
			-
			 M^{\rm x}[C] \le q_{1-\alpha/2} \hat V_n^{\rm x} \}. 
	$$
As  $\sigma^{\rm x} >0$, it follows from Theorem~\ref{thmmain} that
	$$
		\lim_{n \to \infty} \mathbb{P}  \{ M^{\rm x}[C] \in \hat I_n^{\rm x} \} 
	=
		 \mathbb{P}  \{ q_{\alpha/2} \le W \le q_{1-\alpha/2} \}=1 - \alpha. 
	$$
\hfill $\square$

\end{document}